\documentclass[]{siamart171218}

\usepackage{lipsum}
\usepackage{amsfonts}
\usepackage{graphicx}
\usepackage{epstopdf}
\usepackage{algorithmic}
\ifpdf
  \DeclareGraphicsExtensions{.eps,.pdf,.png,.jpg}
\else
  \DeclareGraphicsExtensions{.eps}
\fi

\numberwithin{theorem}{section}

\newcommand{\TheTitle}{Numerical treatment of the nonconservative product 
in a multiscale fluid model for plasmas in thermal nonequilibrium: application to solar physics} 
\newcommand{\TheShortTitle}{Simulation of plasma with nonconservative product} 
\newcommand{\TheAuthors}{Q. Wargnier, S. Faure, B. Graille, T. Magin, M. Massot}

\headers{\TheShortTitle}{\TheAuthors}

\title{{\TheTitle}\thanks{This work was initiated during the 2008 Summer Program at the Center for Turbulence Research at Stanford University. The authors would like to warmly thank Nagi Mansour (NASA Ames Research Center), Fr\'ed\'eric Coquel and Alejandro Alvarez-Laguna (Ecole Polytechnique) for several helpful discussions and for their careful reading of the paper.\funding{The research of Q. Wargnier is funded by an Idex Paris-Saclay interdisciplinary IDI PhD grant, and relies on the the support of NASA Ames Research Center, Advanced Supercomputing Division, von Karman Institute for Fluid Dynamics, EM2C Laboratory at CentraleSup\'elec as well as CMAP and Initiative HPC@Maths from Ecole Polytechnique. T. Magin is supported by a Jean d'Alembert chair of University Paris-Saclay at Ecole Polytechnique.}}}

\author{
  Quentin Wargnier\thanks{CMAP, Ecole polytechnique, CNRS, Universit\'e Paris-Saclay, Route de Saclay,  
91128 Palaiseau Cedex,}
  \and
  Sylvain Faure\thanks{Laboratoire de Math\'ematiques d'Orsay, Univ. Paris-Sud, CNRS, Universit\'e Paris-Saclay, 91405 Orsay, France,}
  \and
  Benjamin Graille\footnotemark[3]
  \and
  Thierry Magin\thanks{von Karman Institute for Fluid Dynamics, Chauss\'ee de Waterloo 72, 1640, Rhode-Saint-Gen\`ese, Belgium.}
	\and
  Marc Massot\footnotemark[2]
 }

\usepackage{amsopn}

\ifpdf
\hypersetup{
  pdftitle={\TheTitle},
  pdfauthor={\TheAuthors}
}
\fi

\usepackage{mathtools}
\usepackage{eufrak}
\usepackage{epstopdf}

\usepackage{graphicx}
\usepackage{caption}
\usepackage{tikz}
\usetikzlibrary{patterns}
\usepackage{pgfplots}
\usepackage{xspace}
\usepackage{xparse} 

\usepackage{graphicx,epstopdf}
\usepackage[caption=false]{subfig}

\usepackage{commandes_GMM}

\begin{document}

\maketitle

\begin{abstract}
This contribution deals with the modeling of collisional multicomponent magnetized plasmas in thermal and chemical nonequilibrium aiming at simulating and predicting magnetic reconnections in the chromosphere of the sun. We focus on the numerical simulation of a simplified fluid model  to investigate the influence on shock solutions of a nonconservative product  present in the electron energy equation. Then, we derive jump conditions based on traveling wave solutions and propose an original numerical treatment in order to avoid non-physical shocks for the  solution that remains valid in the case of  coarse-resolution simulations. A key element for the numerical scheme proposed is the presence of diffusion in the electron variables, consistent with the physically-sound scaling used in the model developed by Graille \textit{et al.}  following a multiscale Chapman-Enskog expansion method [M3AS, 19 (2009) 527--599]. The numerical strategy is assessed in the framework of a solar physics  test case. The computational method is able to capture the traveling wave solutions in both the highly- and coarsely-resolved cases.

\end{abstract}



\section{Introduction} 
\label{sec:sec1}

Plasmas are composed of electrons and heavy particles, such as atoms and molecules,
neutral or ionized. At the microscopic level, electrons and heavy particles do not effectively exchange energy
during collisions, due to their mass disparity. At the macroscopic level, their respective populations of translational energy can be
distributed at different temperatures. Thermal nonequilibrium is found in a variety of plasma applications ranging from astrophysics \cite{alvarez} through electric propulsion \cite{barral} to combustion \cite{tholin}, as well as atmospheric entry flows \cite{magin06}. This contribution deals with the modeling of collisional multicomponent magnetized plasmas in thermal and chemical nonequilibrium, aiming at simulating and predicting magnetic reconnections in the chromosphere of the sun.
If multicomponent magnetohydrodynamic (MHD) simulations are still scarce in solar physics, in recent years, the study of partially ionized plasmas has become an important topic because solar structures, such as prominences \cite{prominences} as well as layers of the solar atmosphere (photosphere and chromosphere) \cite{wiegelmann,Khomenko}, are made of partially ionized plasmas. Multicomponent plasmas introduce physical effects which cannot be described by means of models for fully ionized mixtures, $i.e.$, Cowling's resistivity, isotropic thermal conduction, heating due to ion-neutral friction, heat transfer due to collisions, charge exchange, and ionization, these effects are important to  better understand the behavior of plasmas in the sun chromosphere. Multifluid MHD models \cite{khomenko2,Khomenko,Braginskii} allows us to describe  nonequilibrium effects based on distinct continuity, momentum, and energy equations per species considered in the plasma mixture. However, these models have their own limitations. For instance, the model of Braginskii \cite{Braginskii} is only valid  for high-temperature fully-ionized plasmas. Besides, multifluid models can lead to very stiff systems with characteristic velocities ranging from the speed of sound of the various fluids up to the speed of light \cite{laguna}.  

In this context, Graille \textit{et al.}  \cite{Bible} have derived   from kinetic theory a  multicompo\-nent fluid model  for multicomponent plasmas in thermal nonequilibrium accounting for the influence of the electro-magnetic field, that can be applied  to the sun chromosphere conditions. The model is  obtained  by seeking a generalized Chapman-Enskog solution  to the Boltzmann equation using a multiscale perturbation method. 
A noteworthy difference  with the conventional multifluid models is the scaling  based on a dimensional analysis. This  leads to one single momentum conservation equation and multicomponent diffusion species continuity equations coupled to two equations for the translational (thermal) energy of the electrons and heavy particles. These developments provide a model with an extended range of validity from partially-  to fully-ionized plasmas,  for the non-, weakly-, and strongly- magnetized regimes, together with an entropy inequality and Onsager reciprocity relations for the transport properties. At the  zeroth order of the expansion, this development yields a hyperbolic system of equations with a parabolic regularization of the electron variables due to  dissipative terms such as the electron diffusion velocity and heat flux. When weak solutions (shocks) are considered, the total energy equation is used regarding its conservative form suitable for the development of a numerical scheme. Both the electron and  heavy-particle energy equations exhibit non conservative terms leading to some numerical difficulties that will be reviewed in the next paragraph. In this paper, we propose to select the electron energy equation to close the system, allowing us  to benefit from the regularization of the electron variables. Although our model directly inspired from \cite{Bible} is  directly applicable to the sun chromosphere conditions, it is still necessary to understand how to  treat the nonconservative term present in the hyperbolic convection part of the system. 

Indeed, solving nonconservative hyperbolic systems is a delicate problem because of the definition of weak admissible solutions. 
First, from a theoretical point of view, Dal Maso \textit{et al.} have proposed  in \cite{lefloch2} a theory to define nonconservative products based on the introduction of paths, that generalizes the notion of weak solution for conservative systems  in the sense of distributions. 
In this context, Par\`es  \cite{pares} have developed path-conservative schemes for nonconservative hyperbolic systems. However, it has been proved by Abgrall and Karni \cite{abgrall} that these numerical schemes fail to converge to the right solutions. In fact, even if the correct path is known, the numerical solution obtained can be far from the expected solution, depending mainly on the numerical dissipation of the scheme. In \cite{chalons},  Chalons  and Coquel have proposed  a different strategy for nonconservative hyperbolic systems using Roe-type conservative schemes. They changed the common path-conservative schemes by introducing modified cells in order to compute correctly the solution. Even if progress has been made in the field, the design of accurate and efficient schemes for shock solutions to nonconservative systems of equations and their numerical analysis
still lacks completeness.  Second, from an application point of view, several fields  have encountered this difficulty. 
For two-phase flows, Pelanti and Shyue \cite{pelanti} have proposed an alternative strategy: a Roe solver is used in order to simulate liquid-gas flows with cavitation,  neglecting the nonconservative part of the system. Raviart and Sainsaulieu \cite{sainsaulieu95} have succeeded in evaluating jump conditions  relying on the fact that, in two-phase flows, the nonconservative product acts on linearly degenerate fields. 
In the field of plasma physics, 
Coquel and Marmignon \cite{coquel}  have replaced the equation of  thermal energy of electrons by an equation of conservation of entropy for a model applicable to weakly ionized hypersonic flows in thermal non-equilibrium. Candler and MacCormack  \cite{candler} have considered the nonconservative product in the equation of thermal energy for electrons as a source term for a model applicable to weakly ionized flows. These methods lead to conservative system of equations where the structure of the shock waves is identified, but the link with the original system of equations is still incomplete.
More recently, in \cite{driollet,brull2017}, several numerical schemes have been proposed for the approximation of a nonconservative  compressible Euler system applied to the modeling of  fully ionized plasmas in thermal non-equilibrium; the question of how to evaluate the proper physical jump conditions is not solved.
 In \cite{lowrie}, \cite{edwards}, Lowrie \textit{et al.} have defined semi-analytic solutions for planar radiative shock waves, which can be used for verifying codes for thermal equilibrium diffusion-radiation models. These ideas are also used by Masser \textit{et al.} \cite{masser} to analyze the structure of shock waves in a two-temperature model for fully ionized plasmas. The corresponding ordinary differential equations are integrated and the missing jump relation is obtained by replacing the equation of  thermal energy of electrons by a conservative equation of entropy as in \cite{coquel}, thus avoiding the proper definition and evaluation of a jump condition in the presence of a nonconservative term. This is more difficult in this case, since the nonconservative product acts directly on the genuinely nonlinear waves.
Besides, Shafranov  \cite{shafranov},  Zel\textsc{\char13}dovich and Raizer \cite{zeldovich}, and   Mihalas and Mihalas \cite{mihalas} have been the shock wave structure   for a nonequilibrium fully ionized plasma  in the context of a two-fluid model without nonconservative product.
Relying on the high thermal conductivity of the electrons compared to the one of the ions, the structure of the wave is studied and  the temperature of the electrons is shown to be smooth whereas the temperature of the ions exhibits a discontinuity. In \cite{zeldovich}, it has been shown that the dissipative processes play a major role in the jump conditions for the shock wave: it depends on both the gradients of macroscopic quantities and the transport coefficients. Even if we know what to expect in terms of physics, such waves structure and jump conditions have not been obtained in the framework of a one-fluid model exhibiting a nonconservative product. In summary,  no fully satisfactory solution has yet been achieved to handle both theoretically and numerically the nonconservative product appearing in  a one-fluid model  and it remains an open problem.

In this contribution, we focus on a specific class of solutions written as traveling waves. These solutions are regular enough to ensure that all the terms appearing in the systems are well defined (in particular the nonconservative product). Obviously, the gain of regularity is allowed thanks to the diffusion terms on the electron fields. 
The main goal  is to build a numerical scheme able to capture first the regular solutions of the system and second these particular traveling waves. This is a first step towards capturing the solution of the Riemann problems corresponding to an intermediate scale where the electronic diffusion terms remain. 
The question of the existence of weak solutions (with or without the diffusion terms) is not the subject of this paper. Moreover, we do not deal with the mathematically interesting question of the limit of these traveling waves when the diffusion terms fully vanish since it is beyond the scope of the present paper and it is not relevant for the applications we aim at, where dissipation for the electronic variable always occurs at the same order as the convection phenomena.

We focus on the model derived by Graille \textit{et al.} \cite{Bible} at the zeroth order of the Chapman-Enskog expansion. First, we identify a simplified model, which inherits the same difficulty of dealing with nonconservative products and proper shock numerical solution as the original system, but which is tractable mathematically. A decoupling of the governing equations is proposed and we look for piecewise smooth traveling wave solutions to the decoupled problem, as it is coherent with \cite{zeldovich,shafranov}. This analysis leads to a complete analytical solution, as well as an explicit expression of the missing jump relation for the thermal energy of electrons, where the nonconservative term is to be found. The ability of conducting  the full analysis strongly relies on the proper form of the system and on the presence of a regularizing effect in the electron variables at this order of the expansion. For the numerical solution, we first use a standard finite volume Godunov scheme based on a  consistent discretization of the nonconservative product, in order to resolve the traveling wave. When the level of resolution is too coarse, 
some artificial and non-physical additional shock appears in the solution, whereas for fine resolution, the proper and expected traveling wave solution is reproduced. We thus identify the characteristic scales associated with the compatibility conditions related to the analytic solution for the traveling wave in order to define the limit between the coarse and fine resolutions. Surprisingly enough, it is proved that the smallest diffusion scales associated with the electron mass diffusion have to be properly resolved for the traveling wave to be correctly captured by the numerical scheme. A new scheme based on a specific treatment of the nonconservative product is developed to verify the compatibility equations in a discretized sense. It is able to capture the proper traveling wave even in weakly-resolved cases without generating unexpected additional and artificial numerical shocks. 
Although we use a finite volume Godunov method, the numerical treatment of the nonconservative term and the proposed numerical scheme can be generalized to many finite volume methods: numerical experiments with a Lax-Friedrichs scheme and an upwind scheme are in good agreements. 
The proposed numerical strategy is assessed for a traveling wave test case based on the sun chromosphere conditions, for which the characteristic scales are identified and for which the resolution of the finest diffusion scales is out of reach. Using the new scheme combined to a Strang operator splitting technique, we obtain an accurate resolution of
 the test case with two main advantages: we do not have to resolve the smallest diffusion spatial scales in order to capture the proper traveling wave, as expected, and the timestep is not limited by the Fourier stability condition based on the largest diffusion coefficient.
We eventually investigate the structure of the traveling wave for the original coupled system of equations and prove that: 1- the structure is similar to the one of the decoupled problem, which allow us to have a precious insight on the wave structure and jump conditions; even if we have to resort to a numerical resolution of the missing jump condition by solving a system of 
 ordinary differential equations using a Dormand-Prince (DOPRI853) method with dense output \cite{DORMAND,Hairer93}, we can obtain the missing jump conditions for any Mach numbers in the general case, 2- in a regime of Mach numbers close to one, the missing jump condition for the decoupled and coupled problems are very close to one another, thus fully justifying our strategy to focus on a simplified problem.

The paper is organized as follows: in \Cref{sec:sec2} the model derived by Graille \textit{et al.} \cite{Bible} is presented and briefly compared to other models used by the solar physics community. Then, the decoupling of the governing equations is discussed. In \Cref{sec:sec3},  piecewise smooth traveling wave solutions to the decoupled problem are derived, as well as the missing jump condition associated with the equation of thermal  energy of electrons. The  analytical  solution calculated  is then compared to the ones obtained by solving various models found in the literature. In \Cref{sec:sec4}, a 1D finite volume Godunov scheme with a standard discretization of the nonconservative product is developed, as well as a new scheme based on a specific treatment of the nonconservative product. In \Cref{sec:sec5}, a test case based on the sun chromosphere conditions is fully investigated and our numerical strategy assessed. Finally, in \Cref{sec:sec6}, we show how the ideas developed in the decoupled system case can be extended to the general case and highlight the validity of the decoupled approach in a reasonable Mach number range close to one.

\section{Modeling and governing equations}
\label{sec:sec2}

In this section, we present a multiscale and multicomponent one-fluid model  and explain its advantages and differences compared to the conventional models that are used for describing collisional  plasmas in thermal non-equilibrium, such as those found in the sun chromosphere. We identify a simplified system of equations, which reproduces the difficulty due to the nonconservative product encountered in the general model. An approximation of this system will allow us to conduct analytical studies in \Cref{sec:sec3} and also to highlight the proper paths in the theoretical approach.

\subsection{One simplified model for solar physics}

The non-equilibrium model is  developed based on a thorough kinetic theory derivation by Graille \textit{et al.} \cite{Bible,chimie}. It is   a generalized Chapman-Enskog solution to the Boltzmann equation by means of a multiscale perturbation method. To achieve a fluid limit, the Knudsen number is assumed to scale as a small parameter $\epsilon$ equal to the square root of the ratio of electron mass to a characteristic heavy-particle mass, that drives thermal non-equilbrium. The model is general (see \cref{sec:annexe}) and can be used for multicomponent, non- to weakly- and strongly-ionized plasmas including reactive  collisions between species. 
	
Considering the level of complexity of the general model shown in \cref{sec:annexe} that we propose to use for solar physics applications,  some simplified model is now introduced. We consider the system of equations associated to the zeroth-order of the Chapman-Enskog expansion. At this order, only the electrons have dissipative effects. For simplicity reason, neither Soret-Dufour effects nor electromagnetic forces have been considered. Thermal energy relaxation and chemical processes are assumed to be negligible.  Both heavy particles (multiple species can be considered) and electrons have a common  adiabatic coefficient $\gamma=5/3$, since the internal energy modes are neglected. While diffusion can be anisotropic in the strongly magnetized case, we assume isotropic diffusion since no magnetic field is present. The diffusion structure is still nonlinear even if the electron diffusion coefficient and thermal conductivity are assumed to be constant. 
Under these assumptions, the simplified system of equations  is made up of the conservation equations for the heavy-particle mass, mixture momentum and  total energy, and the electron mass and thermal energy. In nondimensional form (assuming a unit reference Mach number, see \cite{Bible}), this system reads:
\begin{equation}\tag{\ensuremath{S}}
\left\lbrace
\begin{aligned}
&\dt \rhoH +\dx\dscal(\rhoH\vitesse) = 0,   \\
&\dt (\rhoH\vitesse) + \dx\dscal(\rhoH\vitesse\ptens\vitesse + \pression\identite) = 0, \\
&\dt \energie + \dx \dscal(\energie\vitesse +\pression\vitesse) =  \dx\dscal\bigl(\lambda\dx\tempe+\tfrac{\gamma}{\gamma{-}1} \symbolDkl\dx\pre\bigr), \\
&\dt \rhoe + \dx\dscal(\rhoe \vitesse) = \dx\dscal\bigl(\symbolDkl\tfrac{1}{\tempe} \dx \pre\bigr), \\
&\dt (\rhoe\energiee) + \dx\dscal(\rhoe \energiee \vitesse)= 
- \pre\dx\dscal\vitesse + \dx\dscal\bigl(\lambda\dx\tempe+\tfrac{\gamma}{\gamma{-}1} \symbolDkl\dx\pre\bigr) ,
\end{aligned}	
\right.
\label{eq:fullprob}
\end{equation}
where quantity $\rhoH$ stands for the heavy-particle density, $\vitesse$ the heavy-particle velocity, 
$\pression$ the mixture pressure,  $\energie$ the mixture total energy,
$\lambda$ the electron thermal conductivity, $\tempe$ the electron temperature,
$\symbolDkl$, the electron diffusion coefficient, $\pre$ the electron pressure, $\rhoe$, the electron density, $\rhoe\energiee$, the electron thermal energy.
Notice that the mixture pressure is composed of both the electron and heavy-particle partial pressures $\pression=\pre+\pri[\heavy]$, obeying the perfect gas law
\( \pre = (\gamma{-}1)\rhoe\energiee\), \(\pri[\heavy] = (\gamma{-}1)\rhoH\energiei[\heavy]\),
where $\rhoH\energiei[\heavy]$ is the heavy-particle thermal energy. The mixture total energy $\energie$ is defined as
\(
\energie
=
\rhoH\normev^2/2+\pression/(\gamma{-}1)
\).
The term $- \pre\dx\dscal\vitesse$ in the electron energy equation is a nonconservative product examined in great details in this paper. The simplified model \cref{eq:fullprob} inherits the difficulty of dealing with nonconservative products and proper shock numerical solution from the original system, but it remains  mathematically tractable. 

\subsection{Structure of the system}

System~\cref{eq:fullprob} is hyperbolic in the open set of admissible states $\Omega=\{\rhoH>0,\quad\rhoe>0,\quad \vitesse\in \mathbb{R}^3,\quad\pression>0,\quad\pre>0 \}$ with a parabolic regularization on the electron variables. 
For any direction defined by the unit vector $\normale$, the matrix $\normale \dscal  \mathbf{A}$, where   $\mathbf{A}$ is the Jacobian matrix of the  hyperbolic part (removing the second order diffusion terms) is  shown to be diagonalizable with real eigenvalues and a complete set of eigenvectors. The eigenvalue   $\vitesse\dscal \normale$ is  linearly degenerate of multiplicity $\dimension+2$, where $\dimension$ is the space dimension, and the eigenvalues $\vitesse\dscal \normale\pm\soundspeed$ are genuinely nonlinear,  where 
$\soundspeed$ is the sound speed defined by
\(
\soundspeed=\smash{ \sqrt{\gamma\smash{\pression}/\smash{\rhoH}}}
\).
Electrons participate in the momentum balance through the pressure gradient. As a result, the sound speed  includes as well the electron contribution to the pressure.

We recall that  the equation of electron thermal energy is nonconservative in system~\cref{eq:fullprob}. This leads to the difficulties mentioned in the introduction when looking for  discontinuous solutions to the hyperbolic part of the problem. For shock wave solutions, one possibility would be to transform the system~\cref{eq:fullprob} into a conservative system. For instance, the equation of electronic thermal energy can be exchanged with a conservative equation for the electron  entropy \cite{coquel}. This method works only for smooth solutions when there is no dissipation in the electron energy equation. Another possibility would be to consider the nonconservative product as a source term \cite{candler}. However, this strategy modifies the eigenstructure of the system and, as a consequence, the electronic temperature remains constant through a shock wave. Our approach is different: we want to make use of the sound structure of system~\cref{eq:fullprob} in order to derive general jump conditions involving neither simplifications nor modifications of the system. 

\subsection{Approximate decoupled problem}

This section is devoted to an approximate decoupled system obtained by removing the electronic diffusion in the conservation equation for the total energy. This modified system is so-called decoupled since the first three conservation equations constitute the Euler system
and the electronic equations are then solved, once the heavy part velocity $\vitesse$ is known. This simplification has no physical justification since the structure of the diffusion is modified. Nevertheless the decoupled system allows us to derive analytical expressions for traveling wave solutions. The wave structure of the decoupled problem obtained will be shown to be very close to the fully coupled problem of system~\cref{eq:fullprob} in a Mach number regime close to 1, shedding some light on the structure of traveling wave solutions for the fully coupled system. Moreover, this approach allows us to build a numerical scheme which is able to capture the associated traveling wave solutions. Eventually, relying on the same strategy, we will be able to determine the jump conditions for the original five-equation model semi-analytically (using the integration of a dynamical system coupled to a shooting method) and we will check a posteriori that the simplification is physically fully justified in a range of Mach numbers close to one.

The approximate decoupled system splits then into the Euler system
\begin{equation}
\tag{\ensuremath{\montilde{S}_{1}}}
\left\lbrace
\begin{aligned}
&\dt (\wrhoH) +\dx\dscal(\wrhoH\wvitesse) = 0,   \\
&\dt (\wrhoH\wvitesse) + \dx\dscal(\wrhoH\wvitesse\ptens\wvitesse + \wpression\identite) = 0, \\
&\dt \wenergie + \dx \dscal(\wenergie\wvitesse +\wpression\wvitesse) = 0, 
\end{aligned}	
\right.
\label{eq:couple1}
\end{equation}
and a nonconservative drift-diffusion system
\begin{equation}
\tag{\ensuremath{\montilde{S}_{2}}}
\left\lbrace
\begin{aligned}
&\dt\wrhoe + \dx\dscal(\wrhoe \wvitesse) = \dx\dscal\bigl(\tfrac{1}{\wtempe} \symbolDkl\dx \wpre\bigr), \\
&\dt(\wrhoe \wenergiee) + \dx\dscal(\wrhoe \wenergiee \wvitesse)= - \wpre\dx\dscal\wvitesse + \dx\dscal\bigl(\lambda\dx\wtempe+\tfrac{\gamma}{\gamma{-}1}\symbolDkl \dx\wpre \bigr).
\end{aligned}	
\right.
\label{eq:couple2}
\end{equation}
The mixture pressure is introduced as before, $\wpression=\wpre +\wprh$, with the partial pressures 
\(\wpre = (\gamma{-}1)\wrhoe\wenergiee\),
\(\wprh = (\gamma{-}1)\wrhoH\wenergieh\).
We also have $\wpre=\wrhoe\wtempe$.
The mixture total energy is given by
\(\wenergie  = \wrhoH\normewv^2/2+\wpression/(\gamma{-}1)\).
One can find a global solution to this decoupled problem. The first part \cref{eq:couple1} admits discontinuous solutions for which the discontinuity is propagating with the velocity prescribed by the usual Rankine-Hugoniot jump relations and Lax's condition. One can then calculate a solution to the sub-system of electron \cref{eq:couple2} as one piecewise smooth traveling wave, and determine the missing jump condition, with a heavy-particle velocity field previously solved from system \cref{eq:couple1}.

\section{Jump relations and traveling wave solutions to the approximate decoupled problem}
\label{sec:sec3}

In this section, we determine traveling wave solutions for the system~\cref{eq:couple2},  being given the velocity $\wvitesse$ as a piecewise constant function. The variables $(\wrhoH, \wrhoH\wvitesse,\wenergie)$ are assumed to be a shock wave solution to the Euler System~\cref{eq:couple1}
with  velocity  $\sigma$, which  satisfies the Rankine-Hugoniot jump relations and the Lax entropy condition. 
We are subsequently looking for a piecewise smooth traveling wave solution  
in the variables $(\wrhoe,\wenergiee)$, moving with the same velocity $\sigma$, solution to the system~\cref{eq:couple2}. Since the electron variables experience  nonlinear  heat and mass diffusion, their profile is expected to exhibit only weak discontinuities, that is discontinuities in their gradients. We derive boundary conditions at left and right infinities  and show that they do not depend on the (constant) diffusion coefficients $\lambda$ and $\symbolDkl$, hence they can be used as jump conditions associated with the electronic variables $\wrhoe$ and $\wrhoe\wenergiee$. It is consistent with the work of Zel\textsc{\char13}dovich and Raizer in \cite{zeldovich}. These jump conditions are then compared with  literature results. A onedimensional case ($\dimension$=1) is considered.

\subsection{Structure of the traveling wave and jump conditions}

We consider that the heavy-particle variables read as piecewise constant functions depending only on \(\xi = \x\dscal\nunit-\sigma\temps\) 
where $\sigma>0$ is the velocity of the traveling wave prescribed by the Rankine-Hugoniot jump conditions on the heavy-particle variables, and $\nunit$, a unit vector in the first direction, such that $\x\dscal\nunit=x$. 
The same notation is used for functions depending on time and space and for functions depending on $\xi$ as there is no ambiguity. Superscript $\droit$ denotes the right state and $\gauche$ denotes the left state.  We have
\begin{equation}\label{eq:heavyleftrightstates}
\wrhoH(\xi) = \left\lbrace
\begin{aligned}
& \wrhoHL &\text{if }\xi<0,\\
& \wrhoHR&\text{if }\xi>0,
\end{aligned}
\right.
\quad
\wvitesse(\xi) = \left\lbrace
\begin{aligned}
& \wvitesseL &\text{if }\xi<0,\\
& \wvitesseR&\text{if }\xi>0,
\end{aligned}
\right.
\quad
\wpression(\xi) = \left\lbrace
\begin{aligned}
& \wpression^\gauche &\text{if }\xi<0,\\
& \wpression^\droit&\text{if }\xi>0.
\end{aligned}
\right.
\end{equation}
We consider a 3-shock wave that propagates at velocity $\sigma>0$. 
The case of a 1-shock wave is symmetric to the 3-shock case and can be solved in a similar way.

Based on the sub-system of electrons \cref{eq:couple2}, we look for piecewise smooth traveling wave for the electron variables that propagate at the velocity \(\sigma\).
More precisely, the functions \(\wrhoe\) and \(\wpre\) solutions to \cref{eq:couple2} are assumed to satisfy
\begin{itemize}
 \item \(\wrhoe:\xi\mapsto\wrhoe(\xi)\in\Czero(\R)\), \(\wpre:\xi\mapsto\wpre(\xi)\in\Czero(\R)\),
 \item \(\wrhoe\) and \(\wpre\) are \(\Cinf\) on \((-\infty,0)\) and \((0,+\infty)\),
 \item \(\wrhoe\) and \(\wpre\) admit limits in \(\pm\infty\) denoted by
\begin{align*}
&\lim_{\xi \rightarrow -\infty} \wrhoe(\xi)=\wrhoeL,&
&\lim_{\xi \rightarrow +\infty} \wrhoe(\xi)=\wrhoeR,&
&\lim_{\xi \rightarrow \pm\infty}\wrhoe[\prime](\xi)=0,\\
&\lim_{\xi \rightarrow -\infty}\wpre(\xi)=\wpreL,&
&\lim_{\xi \rightarrow +\infty}\wpre(\xi)=\wpreR, &
&\lim_{\xi \rightarrow \pm\infty}\wpre[\prime](\xi)=0.
\end{align*}
\end{itemize}

The goal of this section is to exhibit the structure of these solutions in order to understand the relations between the left and right states according to the traveling wave velocity \(\sigma\) and structure of the diffusion (in particular the value of the coefficients \(\symbolDkl\) and \(\lambda\)).
We assume that the right state $\droit$ is known, and we look for state $\gauche$ connected to the state $\droit$.  

System \cref{eq:couple2} of partial differential equations becomes a system of ordinary differential equations:
\begin{equation}
\left\lbrace
\begin{aligned}
&-\sigma\wrhoe[\prime] + (\wrhoe\wvitesse)' = \symbolDkl (\tfrac{1}{\wtempe}\wpre[\prime])',	\\ 
&-\tfrac{1}{\gamma{-}1}\sigma ~\wpre[\prime] + \tfrac{1}{\gamma{-}1}(\wpre \wvitesse)'= - \wpre \wvitesse' + \lambda~ \wtempe'' + \tfrac{\gamma}{\gamma{-}1}\symbolDkl ~\wpre[\prime\prime] ,
\end{aligned}
\right.	
\label{eq:trav}
\end{equation}
where $\wvitesse$ is a piecewise constant function given in \Cref{eq:heavyleftrightstates}.
The Mach number $\MR$ at state $\droit$ is introduced as 
\(\MR= (\sigma-\wvitesse^\droit)/\cR\),
where $\cR$ is the speed of sound at the right state defined by
\(\cR=\smash{\sqrt{\gamma \smash{\wpression^\droit}/\smash{\wrhoHR}}}\),
and is written in such a way that $\MR>1$ from the Lax conditions. The system \cref{eq:trav} can be solved by considering the two domains $\xi>0$ and $\xi<0$ and their interface $\xi=0$. 
For $\xi>0$, after some algebra, it reads
\begin{equation}
\begin{pmatrix} 
\wpre-\wpreR \\ \wtempe-\wtempe^\droit 
\end{pmatrix}' 
= 
\eta^\droit\begin{pmatrix} 1&-\wrhoeR \\ 
-r^\droit\frac{\gamma-1}{\gamma\wrhoeR}   & r^\droit 
\end{pmatrix}
\begin{pmatrix} 
\wpre-\wpreR \\ 
\wtempe-\wtempe^\droit 
\end{pmatrix},
\label{eq:trav2}
\end{equation}
where the thermal diffusivity $\kappa^\droit$ at the right state and the coefficients $\eta^\droit$ and $r^\droit$ are defined as
\(\kappa^{\droit}=(\gamma-1) \lambda /(\gamma \wrhoeR)\),
\(\eta^\droit=-\cR\MR/\symbolDkl\),
and 
\(r^\droit=\symbolDkl/\kappa^{\droit}\).
Considering that quantities $\kappa^\droit$ and $r^\droit$ are non negative and quantity $\eta^\droit$ non positive for $\gamma>1$, the matrix of \Cref{eq:trav2} has two negative eigenvalues $\delta^{\pm}$ 
\begin{equation}
\delta^{\pm}=\tfrac{1}{2}\eta^\droit \Bigl(1+r^\droit \pm\sqrt{(1+r^\droit)^2-\tfrac{4}{\gamma}r^\droit}\Bigr).
\label{eq:eigenvalue}
\end{equation}

Finally, for $\xi>0$ one gets an analytical expression of the solution that combines decreasing exponential functions 
\begin{align*}
\wpre(\xi) &= \wpreR + \wrhoeR\bigl(K^{R+}\e^{\delta^{+}\xi} + K^{R-}\e^{\delta^{-}\xi} \bigr),&&\xi>0,\\
\wtempe(\xi) &= \wtempe^\droit
+ \bigl(1-\tfrac{\delta^+}{\eta^\droit}\bigr)K^{R+}\e^{\delta^{+}\xi}
+ \bigl(1-\tfrac{\delta^-}{\eta^\droit}\bigr)K^{R-}\e^{\delta^-\xi}, &&\xi>0,
\end{align*}
where $K^{R\pm}$ are two integration constants that need to be determined by using the continuity and the jump of the derivative gradients in \(\xi=0\). 
Moreover, for $\xi<0$, similar algebraic relations as those found in \Cref{eq:trav2} are obtained by replacing state $\droit$ by state $\gauche$. It leads to similar eigenvalues as those in \Cref{eq:eigenvalue} by replacing state $\droit$ by state $\gauche$. However, the only way for having a non diverging solution is to get constant functions equal to the left constant state $\gauche$ by setting the integration constants $K^{L\pm}=0$. Indeed, the only bounded solutions when \(\xi\) goes to \(-\infty\) are constant solutions.

At \(\xi=0\), since $\wpre$ and $\wrhoe$ (and $\wtempe$) are continuous functions, the system~\cref{eq:trav} leads to the following compatibility equations
 \begin{equation}
\wpre(0)\saut{\wvitesse} = \symbolDkl\saut{\wpre[\prime]}, 
\qquad
\tfrac{\gamma}{\gamma{-}1} \wpre(0)\saut{\wvitesse} = \lambda\saut{\wtempe'} + \tfrac{\gamma}{\gamma{-}1}\symbolDkl\saut{\wpre[\prime]},
\label{eq:discontinuity}
\end{equation}
where $\saut{\cdot}$
denotes for the value of the jump in $\xi=0$. 
The term of the left hand side in the second equality of \Cref{eq:discontinuity} is the contribution of two terms: a first one coming from the convective part and a second one coming from the nonconservative product. Actually, since $\wpre$ is continuous and the derivative $\wvitesse'$ has a jump at \(\xi=0\), the nonconservative product in the second equation of \Cref{eq:trav} is not ambiguous. Finally, the second relation of \Cref{eq:discontinuity} can be simplified by using the first one. The compatibility \Cref{eq:discontinuity} then become
\begin{equation}\label{eq:discontinuitysimplified}
\wpre(0)\saut{\wvitesse} = \symbolDkl \saut{\wpre[\prime]} ,
\qquad
\saut{\wtempe'}=0.
\end{equation}

The second relation of \Cref{eq:discontinuitysimplified} provides the continuity of the derivative of $\tempe$ in the discontinuity (at \(\xi=0\)). In other words, $\tempe$ is a $\Cun$ function. This result is consistent with the work of Zel\textsc{\char13}dovich and Raizer in \cite{zeldovich}, showing that the temperature of electron is smooth in the shock wave of a nonequilibrium fully ionized plasma.
Moreover, the first relation of \Cref{eq:discontinuitysimplified} can be seen as a relation between the jump of the pressure gradient in the discontinuity, the jump of the velocity and the diffusion coefficient $\symbolDkl$. This result is also consistent with the work of Zel\textsc{\char13}dovich and Raizer in \cite{zeldovich}, showing that the discontinuity of the shock wave in a plasma depend on the dissipative process.
Defining the two characteristic lengths of the diffusion $\LD$ associated to the electron diffusion coefficient and $\LT$ associated to the thermal conductivity by
\begin{equation}
\LD
=
\frac{\symbolDkl}{\saut{\wvitesse}}, 
\qquad 
\LT
=
\frac{\gamma{-}1}{\gamma}
\frac{\lambda}{\wrhoeR}\frac{1}{\left[\wvitesse\right]_{(0)}}=\frac{\kappa^\droit}{\left[\wvitesse\right]_{(0)}},
\qquad
\left[\wvitesse\right]_{(0)}>0,
\label{eq:LD}
\end{equation} 
\Cref{eq:discontinuitysimplified} can be rewritten as
\begin{equation}
\saut{\wpre[\prime]} =  \tfrac{1}{\LD} \wpre(0),
 \qquad
 \saut{\wtempe'} = 0.
\end{equation}

\begin{figure}[htbp]
\center{
\tikzset{
    fleche/.style = {->, thick},
    courbe/.style = {mark=none, very thick, samples=100},
    pente/.style = {mark=none, thick, dashed},
}
\begin{tikzpicture}[scale = .6]
\begin{axis}[
  grid=major,
  xtick={0, 1},
  xticklabels={0, $\LD$},
  ytick={0, .3, 1},
  yticklabels={0, $\wpreR$, $\wpreL$},
  xlabel={$\xi=x-\sigma t$},
  ylabel={$\wpre$},
  xlabel style={below},
  ylabel style={right},
  xmin=-5,
  xmax=5.5,
  ymin=-0.1,
  ymax=1.25,
  ]
\addplot+[courbe, draw=black, domain=-5:0] {%
	1%
	};
\addplot+[courbe, draw=black, domain=0:5] {%
	.3+.7*exp(-x)
	};
\addplot+[pente, draw=black, domain=0:1] {%
	1-x
	};
\end{axis}
\end{tikzpicture}
\hfill
\begin{tikzpicture}[scale = .6]
\begin{axis}[
  grid=major,
  xtick={0, 1},
  xticklabels={0, $\LD$},
  ytick={0, .4, 1},
  yticklabels={0, $\wtempe^\droit$, $\wtempe^\gauche$},
  xlabel={$\xi=x-\sigma t$},
  ylabel={$\wtempe$},
  xlabel style={below},
  ylabel style={right},
  xmin=-5,
  xmax=5.5,
  ymin=-0.1,
  ymax=1.25,
  ]
\addplot+[courbe, draw=black, domain=-5:0] {%
	1%
	};
\addplot+[courbe, draw=black, domain=0:5] {%
	.4+.6*exp(-x*x)
	};
\addplot+[pente, draw=black, domain=0:1] {%
	1-x
	};
\end{axis}
\end{tikzpicture}
}
\caption{Scheme of the traveling wave with the characteristic diffusion length $\LD$}
\label{fig:schema}
\end{figure}
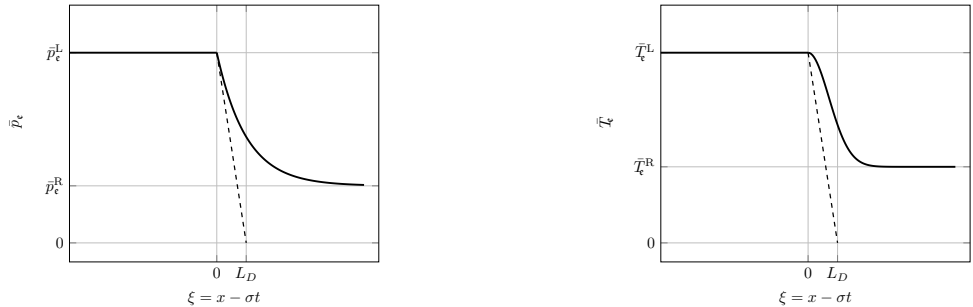

The jump compatibility relations for the pressure of electron $\wpre$ and for the electron temperature $\wtempe$
 link the  states $\droit$ and $\gauche$ at infinity. They are obtained by integrating  \Cref{eq:trav} from $0^+$ to $\infty$, then using \Cref{eq:discontinuitysimplified} and the jump conditions for the heavy particle velocity $\wvitesse$.
Finally, these relations can be combined to obtain the density jump conditions between the  states $\droit$ and $\gauche$  at infinity.
One gets
 \begin{equation}
\frac{\wpreL}{\wpreR}=\frac{(\gamma+1)\MR^2}{(1-\gamma)\MR^2 + 2\gamma}, 
\quad 
\frac{\wtempe^\gauche}{\wtempe^\droit}=\frac{(\gamma-1)\MR^2 + 2}{(1-\gamma)\MR^2+2\gamma},
\quad
 \frac{\wrhoeL}{\wrhoeR}  = \frac{(\gamma+1)\MR^2}{(\gamma-1)\MR^2 + 2}.
\label{eq:jumppe}
\end{equation}
The state $\gauche$ at infinity does not depend on the diffusion coefficients $\symbolDkl$ and $\lambda$. However, according to \Cref{eq:discontinuitysimplified}, the jump compatibility relations depend on the variables and their gradients in the discontinuity. This result is consistent with the work of Zel\textsc{\char13}dovich and Raizer in \cite{zeldovich}, and Shafranov in \cite{shafranov}. 
Notice that the jump condition ${\wrhoeL}/{\wrhoeR} $ is the same as for ${\wrhoHL}/{\wrhoHR}$ and is compatible with the Rankine-Hugoniot jump relations. 
Let us underline however that the relations of \Cref{eq:jumppe} are valid only for $\MR^2 < 2\gamma/(\gamma-1)$, 
and a singularity is present when $\MR$ goes to $\sqrt{2\gamma/(\gamma-1)}$. 
In the case where the Mach number is equal to this value, it is not possible to solve the problem of the traveling wave and build a solution in the expected form. In addition, beyond this limit value, the solution is showing negative temperatures. Thus the decoupling of the system provides a reasonable solution only below this value. In fact this coupling of heavy species and electrons is expected to be of very weak amplitude in the neighborhood of Mach one since the dissipation is going to be really weak for such small amplitude shocks; 
dissipation is only playing a role in that Mach range in order to regularize the electrons profiles and reach the proper jump conditions on the electronic temperature. In particular the heavy particle decoupled hydrodynamic velocity jump is a good prediction in this range of the velocity jump condition for the full system, thus justifying the weak coupling leading to the decoupled system. When the Mach number is increasing, the full system exhibit a strong coupling between the jump of electronic variables and the jump in heavy particle hydrodynamic velocity, which can not be reproduced by the decoupled system and eventually leads to an incompatibility in terms of jump conditions of the electronic variables.
Thus, the obtained relations are valid for a Mach number range close to one and will prove, in Section \ref{sec:sec6}, to provide a good estimate in that range of the jump conditions for the full system.

\subsection{Comparison with classical jump conditions in the literature}

In this section, the jump conditions proposed in \Cref{eq:jumppe} are compared with other usual jump conditions from conservative system of equations.

First, one can consider a conservative system of equations where we replace the nonconservative equation of electron internal energy by an equation of conservation of electron entropy (see \Cref{eq:entropy} in \Cref{sec:annexe2}). The obtained model is called model $\mathcal{M}^{\text{ent}}$. In this case, one would get the following jump conditions:
\begin{equation}
\frac{\wpreL}{\wpreR} 
\underset{\mathcal{M}^{\text{ent}}}{\;:=\;}
\left(\frac{(\gamma+1)\MR^2}{(\gamma-1)\MR^2+2}\right)^{\gamma},
\quad 
\frac{\wtempe^\gauche}{\wtempe^\droit}
\underset{\mathcal{M}^{\text{ent}}}{\;:=\;}
\left(\frac{(\gamma+1)\MR^2}{(\gamma-1)\MR^2+2}\right)^{\gamma-1}.
\label{eq:jumppee}
\end{equation}

Second, one can consider another conservative system of equations where the nonconservative product is considered as a source term: only the conservative part of the system is considered for getting the jump conditions (see \Cref{eq:source} in \Cref{sec:annexe2}). The model obtained in that case is called model $\mathcal{M}^{\text{src}}$ and the jump conditions read
\begin{equation}
\frac{\wpreL}{\wpreR}
\underset{\mathcal{M}^{\text{src}}}{\;:=\;}
\frac{(\gamma+1)\MR^2}{(\gamma-1)\MR^2+2}, 
\quad 
\frac{\wtempe^\gauche}{\wtempe^\droit}
\underset{\mathcal{M}^{\text{src}}}{\;:=\;}
1.
\label{eq:jumppes}
\end{equation}

\begin{figure}[h!]
\centering
\hfill
\includegraphics[scale=0.4]{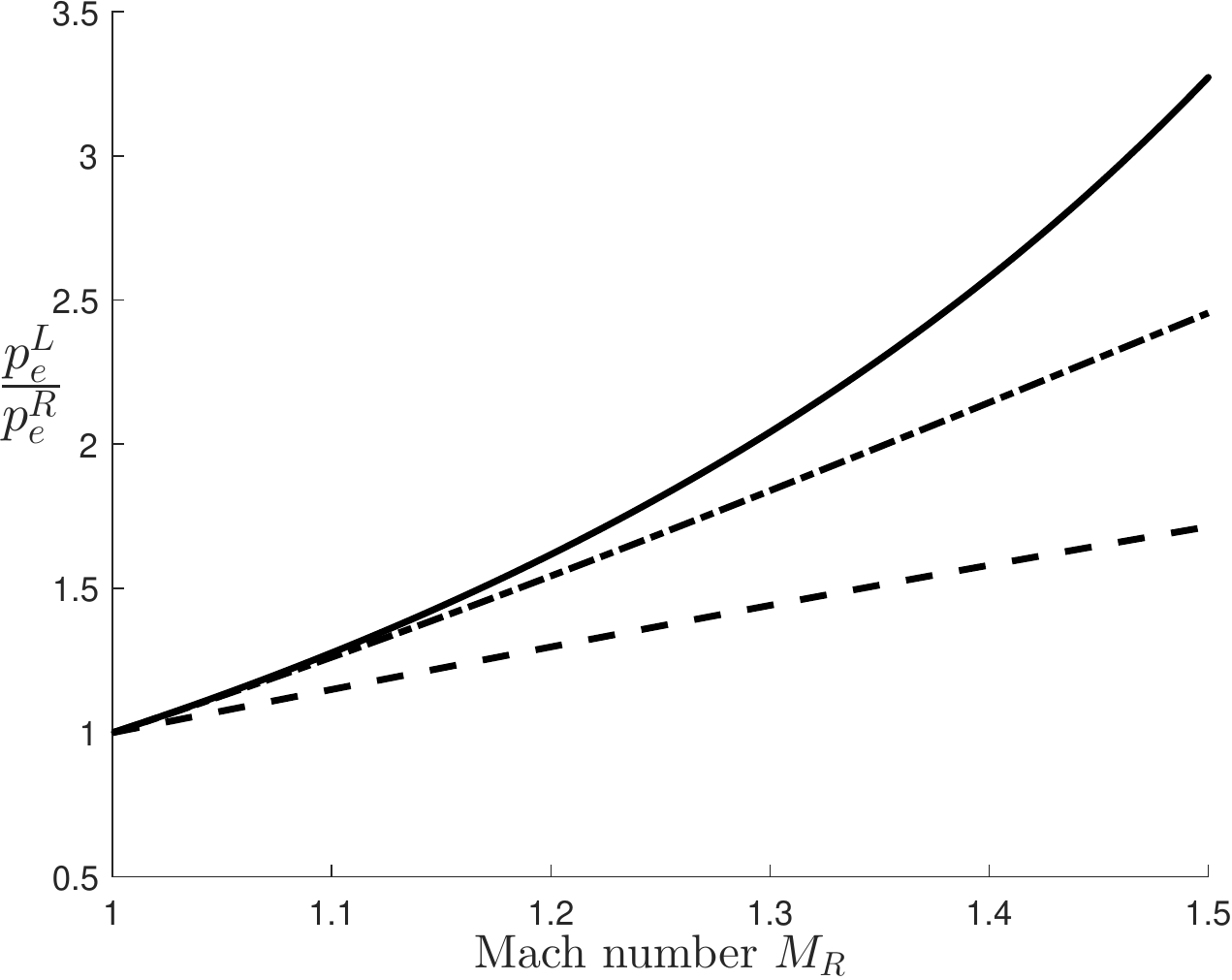}
\hfill
\includegraphics[scale=0.4]{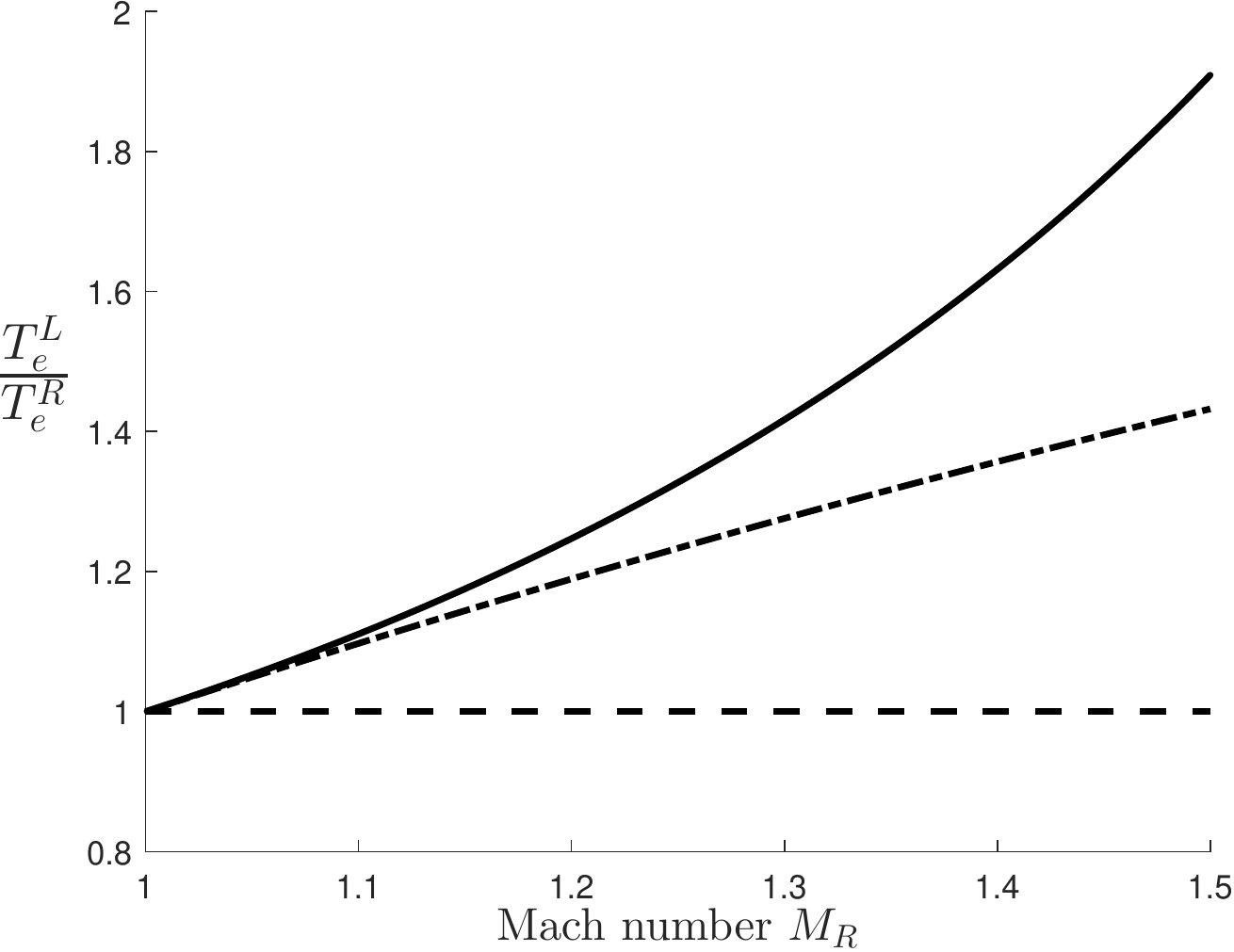}
\hfill \null
\caption{Ratio $\wpreL/\wpreR$ and $\wtempe^\gauche/\wtempe^\droit$ as a function of the Mach number $\MR$, from \Cref{eq:jumppe} in full line, from \Cref{eq:jumppee} in semi-dashed line and from \Cref{eq:jumppes} in dashed line.}
\label{fig:resultsjump}
\end{figure}

The jump conditions for the electron pressure and electron temperature obtained by the traveling wave method in \Cref{eq:jumppe} are then compared to those obtained by means of \Cref{eq:jumppee}  and \Cref{eq:jumppes}.
In \cref{fig:resultsjump}, the three jump conditions are plotted as functions of the Mach number
for  Mach numbers between $1$ and $1.5$.
We observe first that the isothermal jump conditions of model $\mathcal{M}^{\text{src}}$ rapidly underestimates the post-jump temperature. Moreover, this model is not reasonable since the dynamics of smooth waves, such as rarefaction waves,  is modified. 
Second, the jump conditions \Cref{eq:jumppes} of model $\mathcal{M}^{\text{ent}}$ are similar to the ones of \Cref{eq:jumppe} for a Mach number regime close to 1. However, significant differences can be observed when the Mach number is increasing.

As a conclusion, by looking for piecewise smooth traveling wave solutions, we were able to get an analytical expression of the missing jump condition associated with the thermal energy of electrons, and the analytical traveling wave solution of the electron variables valid for a Mach number regime close to one. For that purpose, we had to decouple the problem: a discontinuity propagating at velocity $\sigma$ where the jump conditions are prescribed by the usual Rankine-Hugoniot relations solution to the sub-system \cref{eq:couple1} and a continuous traveling wave propagating at the same velocity $\sigma$, solution to the sub-system \cref{eq:couple2}.  The resulting jump conditions are valid in a neighborhood of Mach one and then lead to a singularity for larger Mach numbers. 
In the interval where they are valid, they exhibit rather important differences with the conditions found in the literature. 
In section \ref{sec:sec6}, we will prove that the analytical expression obtained for the decoupled system is a very good approximation in a Mach range close to the one of the jump conditions for the fully coupled problem \cref{eq:fullprob}, which does not lead to any singular behavior. The jump condition for the fully coupled problem \cref{eq:fullprob} will also be proved to be very different from the usual jump conditions of the literature.  The next step is to verify numerically the jump conditions and if we can capture the analytical traveling waves. 

\section{Numerical scheme for the decoupled system \cref{eq:couple1,eq:couple2}}
\label{sec:sec4}

In the previous section, we show the existence of a traveling wave for system~\cref{eq:couple2}. The aim of this section is to develop numerical methods able to capture the dynamic of the traveling wave. First, we introduce a standard scheme based on a finite volume Godunov method with a standard discretization of the nonconservative product and then a specific treatment. 

Note that the proposed method to treat the nonconservative term is independant of the chosen finite volume scheme. Numerical experiments for a Lax-Friedrichs scheme and an upwind scheme have been performed and led to the same conclusions. 

\subsection{Finite volume scheme with a standard discretization for the nonconservative product}

The electron variables $\wrhoe$ and $\wrhoe\wenergiee$ are initialized by using the analytical solution found in the previous section. The heavy variables $\wrhoH$, $\wvitesse$, and $\wenergie$ are discontinuities propagating at velocity $\sigma$, where the  conditions are fixed by the Rankine-Hugoniot conditions.

\begin{figure}[htbp]
\center{\small
\tikzset{
    delta/.style = {<->},
    fleche/.style = {->},
    trait/.style = {mark=none, very thin},
    choc/.style = {mark=none},
    detente/.style = {mark=none, very thin, fill=gray!25, opacity=0.5},
    discontinuity/.style = {mark=none, dashed},
    centre/.style = {fill=gray!25},
    cellule/.style = {very thick},
}
\begin{tikzpicture}[scale = .7]

\draw[cellule] (-1,0) -- (1,0);
\draw[cellule] (-1,2) -- (1,2);

\draw[trait, fleche] (-3.5,0) -- (3.5,0) node[below]{$\x$};
\draw[trait] (-3.5,2) -- (3.5,2);
\draw[trait, fleche] (-3,0) -- (-3,3) node[right]{$t$};
\draw[trait] (-1,0) -- (-1,2.5);
\draw[trait] (1,0) -- (1,2.5);
\draw[trait] (3,0) -- (3,2.5);

\draw[centre] (-2,0) circle (0.1);
\draw (-2,-.75) node{$\cellj[j-1]$};
\draw (-1,-.25) node{$x_{j-\onehalf}$};
\draw[centre] (0,0) circle (0.1);
\draw (0,-.75) node{$\cellj$};
\draw (1,-.25) node{$x_{j+\onehalf}$};
\draw[centre] (2,0) circle (0.1);
\draw (2,-.75) node{$\cellj[j+1]$};

\draw[trait, delta] (-1,2.75) -- (1,2.75);
\draw (0,3) node{$\Dx$};
\draw[trait, delta] (-4,0) -- (-4,2);
\draw (-4.5,1) node{$\Dt$};

\draw[detente] (-1,0) -- (-2,2.5) -- (-1.25,2.5) -- cycle;
\draw[discontinuity] (-1,0) -- (-.75,2.5);
\draw[choc] (-1,0) -- (0.,2.5);
\draw[choc] (1,0) -- (0.25,2.5);
\draw[discontinuity] (1,0) -- (1.45,2.5);
\draw[choc] (1,0) -- (2.,2.5);

\end{tikzpicture}
}
\caption{Notations for finite volume scheme.}
\label{fig:godunov}
\end{figure}
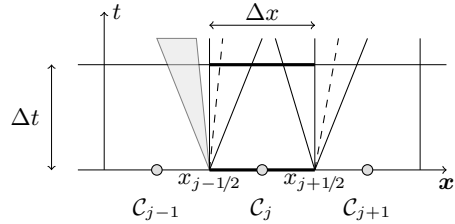

A monodimensional finite volume Godunov method is used to discretize the electron sub-system \cref{eq:couple2}. We consider a finite domain of length $L$ with $N$ cells of length $\Dx = L/N$. The position of each cell 
$\cellj$, $1\leq j\leq N$ is defined by its center $x_j$ at the middle of the interfaces $x_{j+\onehalf}$ and $x_{j-\onehalf}$. 
The bounds of the domain are not taken into account as the simulations are stopped before any interaction occurs between the traveling wave and the boundary. Left and right Dirichlet conditions are then used.
The time is also discretized with a timestep $\Dt$. 
\Cref{fig:godunov} can be used to visualize these standard notations. 

We denote by $\Ujn$, $n\geq0$, $1\leq j\leq N$, the vector of the natural variables at time $t^n$ in the cell $\cellj$ is 
\begin{equation*}
\Ujn=\left( \Ujnun, \Ujndeux \right),\quad
 \Ujnun=\wrhoejn \text{ and } \Ujndeux = \wrhoejn\wenergieejn.
\end{equation*}
The general scheme reads
\begin{equation}\label{eq:generalscheme}
\Ujnp= \Ujn -\tfrac{\Dt}{\Dx} \bigl(\Fjcp - \Fjcm \bigr) + \tfrac{\Dt}{\Dx} \bigl(\FjDp - \FjDm \bigr) + \Fjnc ,
\end{equation}
where 
$\Fjcpm$ are the convective fluxes at interfaces $j\pm\onehalf$, 
$\FjDpm$ the diffusive fluxes at interfaces $j\pm\onehalf$, 
and $\Fjnc$ the value of the nonconservative term in the cell~$j$. 

The convective flux $\Fjcp$ is computed by means of Godunov's scheme, by solving the Riemann problem with the left and the right values given by the cells $j$ and $j{+}1$ and by taking the value of the flux at the interface. 
Note that the solution of the Riemann problem used by this Godunov's solver is essentially the solution of the transport equation with a constant velocity but at the interfaces close to the discontinuity of the velocity $\wvitesse$ (in these cases, the traveling wave is used). 

Then, the diffusive flux is calculated by using a second-order centered scheme 
\begin{equation} \label{eq:diffusiveflux }
 \FjDp 
 =\tfrac{1}{\Dx}
 \Bigl(
 D  \tfrac{\gamma - 1}{\tempepn}
 (\Ujndeux[j+1]-\Ujndeux[j]) 
 ,
 \lambda ( \tempejpn - \tempejn) 
 + 
 D \gamma (\Ujndeux[j+1]-\Ujndeux[j])
 \Bigr),
\end{equation}
where $\tempejn=(\gamma-1)\wenergieejn=(\gamma-1)\Ujndeux/\Ujnun$ and the interface temperature $\tempepn$ reads
\begin{equation*}
 \tempepn = \frac{\gamma-1}{2} \Bigl(
 \frac{\Ujndeux[j+1]}{\Ujnun[j+1]} + \frac{\Ujndeux[j]}{\Ujnun[j]}
 \Bigr) = \frac{1}{2} \Bigl(\tempejpn+\tempejn \Bigr).
\end{equation*}
Note that the error of consistency of this second-order discretization can be bounded by a \(\mathcal{O}(\Dx^2)\) term. This accuracy is sufficient for the proposed scheme as the major error is done by the Godunov's discretization of the hyperbolic part.

Finally, the second composant of the nonconservative term $\Fjnc = (0, \Fjncdeux)$ is computed as an approximation (as accurate as possible) of the integral over $[t^n,t^{n+1}]\times\cellj$ of the nonconservative contribution
\begin{equation*}
 \int_{t^n}^{t^{n+1}}\!\! \int_{\cellj} \wpre \dx\wvitesse \operatorname{d}\!\x  \operatorname{d}\!t  
 \simeq
 (\gamma-1)\Ujndeux
 \int_{t^n}^{t^{n+1}}\!\! \int_{\cellj}  \dx\wvitesse \operatorname{d}\!\x  \operatorname{d}\!t ,
\end{equation*}
where the integral of $\dx\wvitesse$ is exactly computed as the velocity $\wvitesse$ is prescribed.
The proposed scheme is then a consistent numerical scheme with a standard discretization of the nonconservative product that can be tested for capturing the traveling wave.

\subsection{Numerical results}
\label{sec:num}
In this section, we present some numerical experiments, using the finite volume scheme with standard discretization for the non-conservative product proposed in the previous section. Different resolutions of the traveling wave are presented with a double objective: 
first, to capture the dynamic of the traveling wave with a fine enough mesh;
second, to visualize the behaviour with a coarse mesh in order to understand how the scheme can capture a shock discontinuity.
We focus on the 3-wave with respect to the wave structure of the Euler system \cref{eq:couple1}, so that the right state $R$ is known. Besides, a supersonic regime is studied for a Mach number close to one in order to guarantee the existence of the traveling wave introduced in \Cref{sec:sec3}.

The number of nodes $N$ and the length of the domain $L$ are fixed: $N=2000$ and $L=10$. 
The initial position of the traveling wave (that is the position corresponding to $\xi=0$ in the moving frame) is $0.2L$. 
The time discretization $\Dt$ is fixed by a Fourier condition $\Dt\leq\tfrac{1}{2}\beta\Dx^2$ where $\beta=\max(D,\kappa^\droit)$. 
All the simulations are stopped at $t =t_f= 1$, corresponding to a displacement of the traveling wave to $0.373L$. The electron thermal conductivity $\lambda=0.001$ is fixed, so the associated characteristic length $\LT=7.6\times 10^{-2}$ is fixed. This value has been chosen as a good compromise between the length of the domain, which is fixed, and the regularization of the profile of the electron temperature. However, the electron diffusion coefficient $\symbolDkl$ is going to vary in our numerical experiments changing the resolution of the traveling wave. This diffusion length $\LD$ is related to the diffusion coefficient $\symbolDkl$ in \Cref{eq:LD} as an increasing function. 
Consequently, since the length of the domain and the number of nodes are fixed, we improve the resolution of the traveling wave in the characteristic length $\LD$ by increasing the diffusion coefficient $\symbolDkl$. Simulations are conducted for diffusion coefficient $\symbolDkl$ between $10^{-3}$ and $10^{-1}$ corresponding to different resolutions of the traveling wave and lengths $\LD$. The two extreme cases, denoted by case \caseHD (over-resolved)
and case \caseWD (under-resolved),
are reported in \Cref{tab:tablea}: 
in case \caseHD, $\LD> \LT$ whereas in case \caseWD,  $\LD\ll \LT$. In \Cref{sec:sec5}, we will see that the physical test case in sun chromosphere conditions corresponds to cases where $\LD\ll \LT$.

\begin{table}[ht]
\centering\footnotesize
  \caption{Values of the diffusion coefficient $\symbolDkl$ used in the numerical experiments}
  \label{tab:tablea}
\grandtraittop
\renewcommand*{\arraystretch}{1.1}
\begin{tabular}{@{}p{4cm}p{2cm}p{2cm}@{}}
&case \caseHD&case \caseWD \\ \hline
$\symbolDkl$ & $10^{-1}$ & $10^{-3}$ \\
$\LD$ & $3.055{\times}10^{-1}$ & $3.055{\times}10^{-3}$ \\
Number of nodes in $\LD$ &61.1&0.611
\end{tabular}
\grandtraitbottom
\end{table}
	
The right state and the left state have been initialized with the values given in \Cref{tab:tableb}. The left state of the traveling wave has been computed using Rankine-Hugoniot relations for the heavy particles variables and jump conditions given in \cref{eq:jumppe} for the electronic variables.

\begin{table}[ht]
\centering\footnotesize
  \caption{Right and left states of the traveling wave}
  \label{tab:tableb}
\grandtraittop
\renewcommand*{\arraystretch}{1.1}
\begin{tabular}{@{}p{2cm}p{1.1cm}p{1.1cm}p{1.1cm}p{1.1cm}p{1.1cm}p{1.8cm}@{}}
&$\wrhoH$ & $\wrhoe$ & $\wpression$&$\wpre$&$\wvitesse$&Mach number \\ \hline
right state $R$ & 1 &0.01 &1&0.1&0.2&1.1832 \\
left state $L$ &1.274&0.01274&1.5&0.1556&0.527&0.8563
\end{tabular}
\grandtraitbottom
\end{table}

%

\begin{figure}[tbhp]
\centering
\subfloat[Case \caseHD]{\label{fig:figcase1}\includegraphics[scale=0.25]{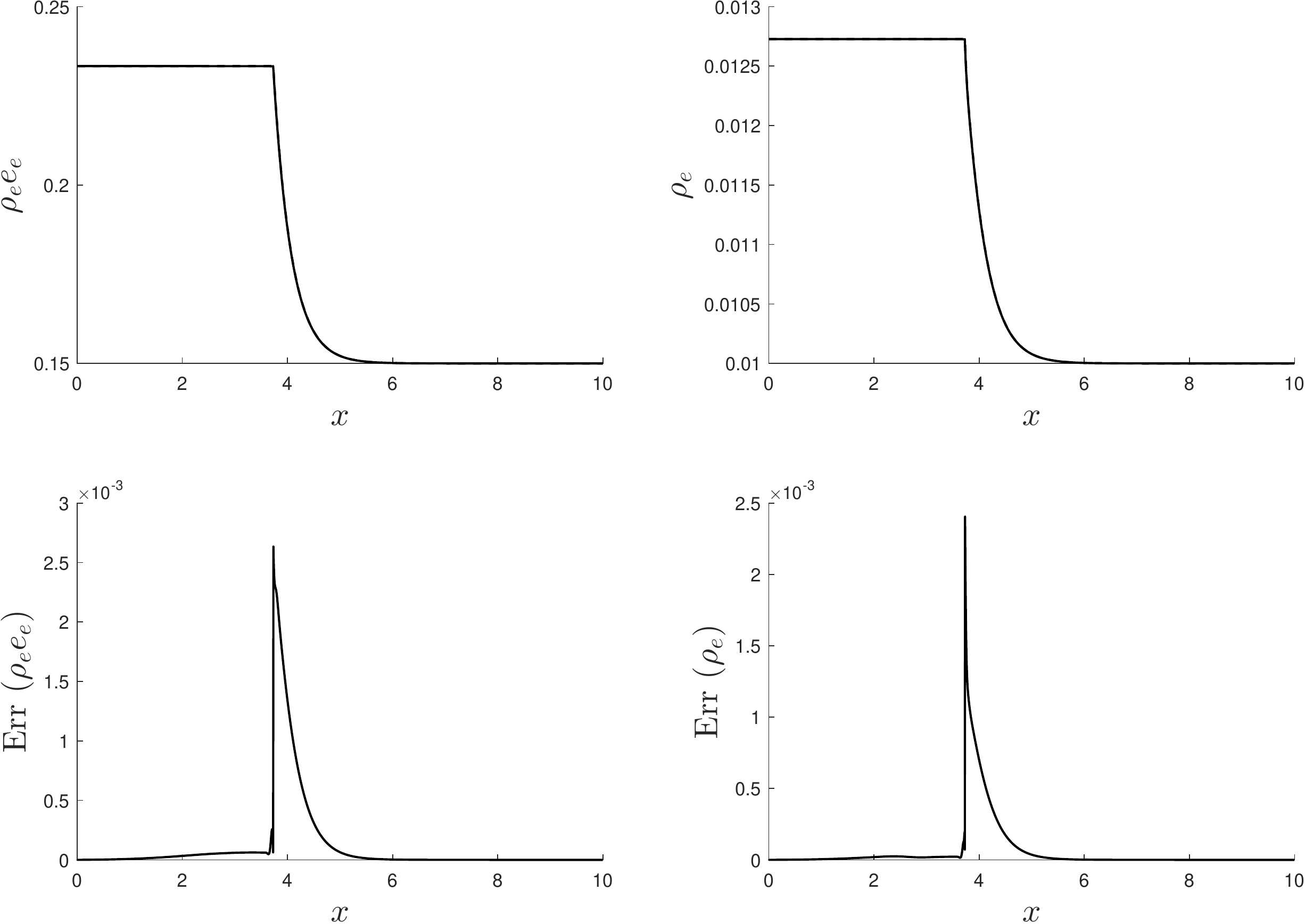}}
\hfill
\subfloat[Case \caseWD]{\label{fig:figcase2}\includegraphics[scale=0.25]{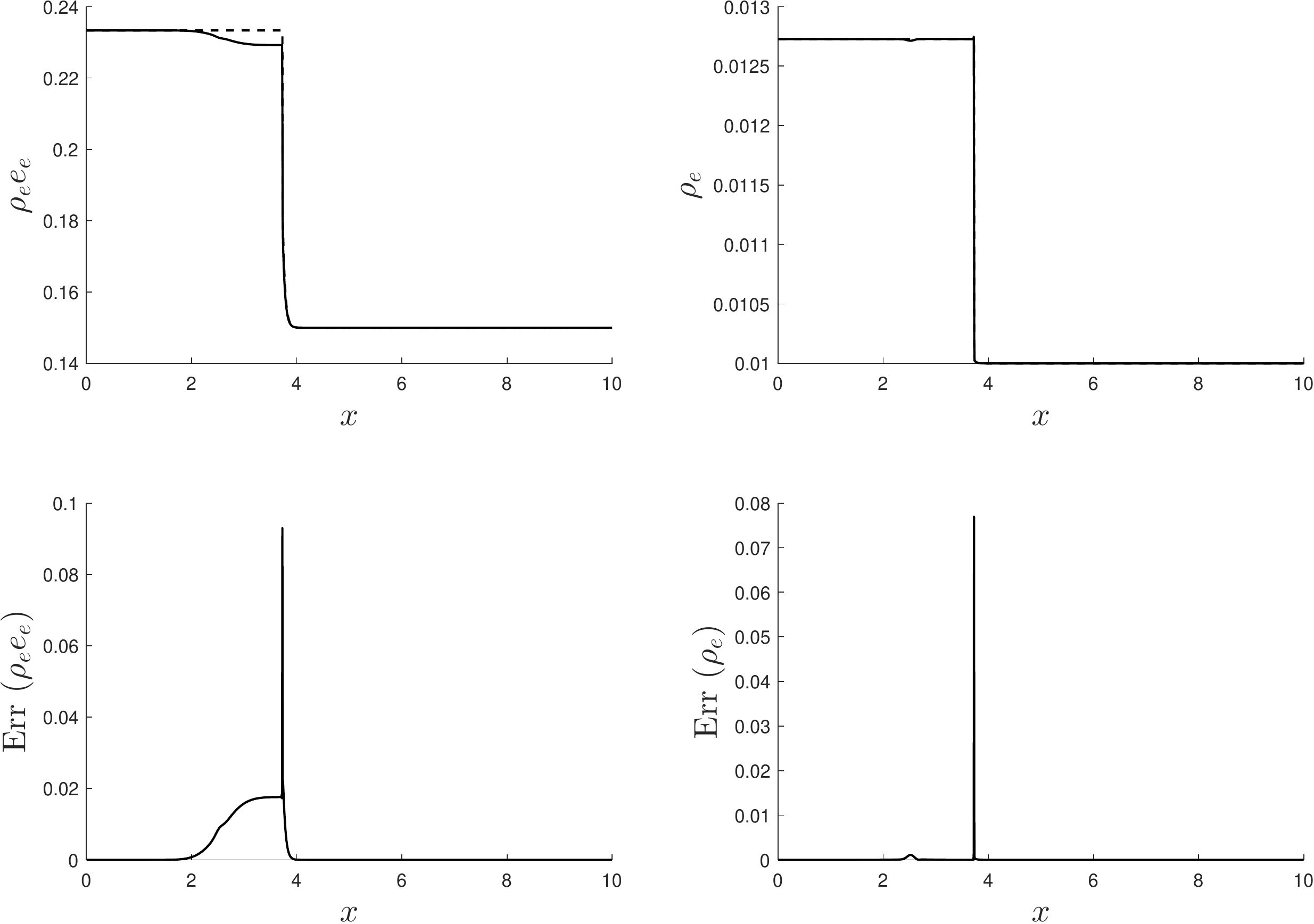}}
\caption{Analytical solution (dashed line), numerical solution (full line), and relative error (Err) for the traveling wave   at $t=t_f=1$}
\end{figure}			

\Cref{fig:figcase1,fig:figcase2} show first that the traveling wave is well captured in the case~\caseHD: the dynamic of the traveling wave is preserved if the number of nodes in the length $\LD$ is large enough; second, that a non-expected artificial numerical shock appears in the case~\caseWD. Two main contributions to the difference between the numerical and analytical solutions can be exhibited:
a contribution upstream of the shock due to the numerical dissipation in the regular part of the traveling wave;
a contribution downstream of the shock due to the error on the gradients in the discontinuity.

\begin{figure}[tbhp]
\includegraphics[scale=0.4]{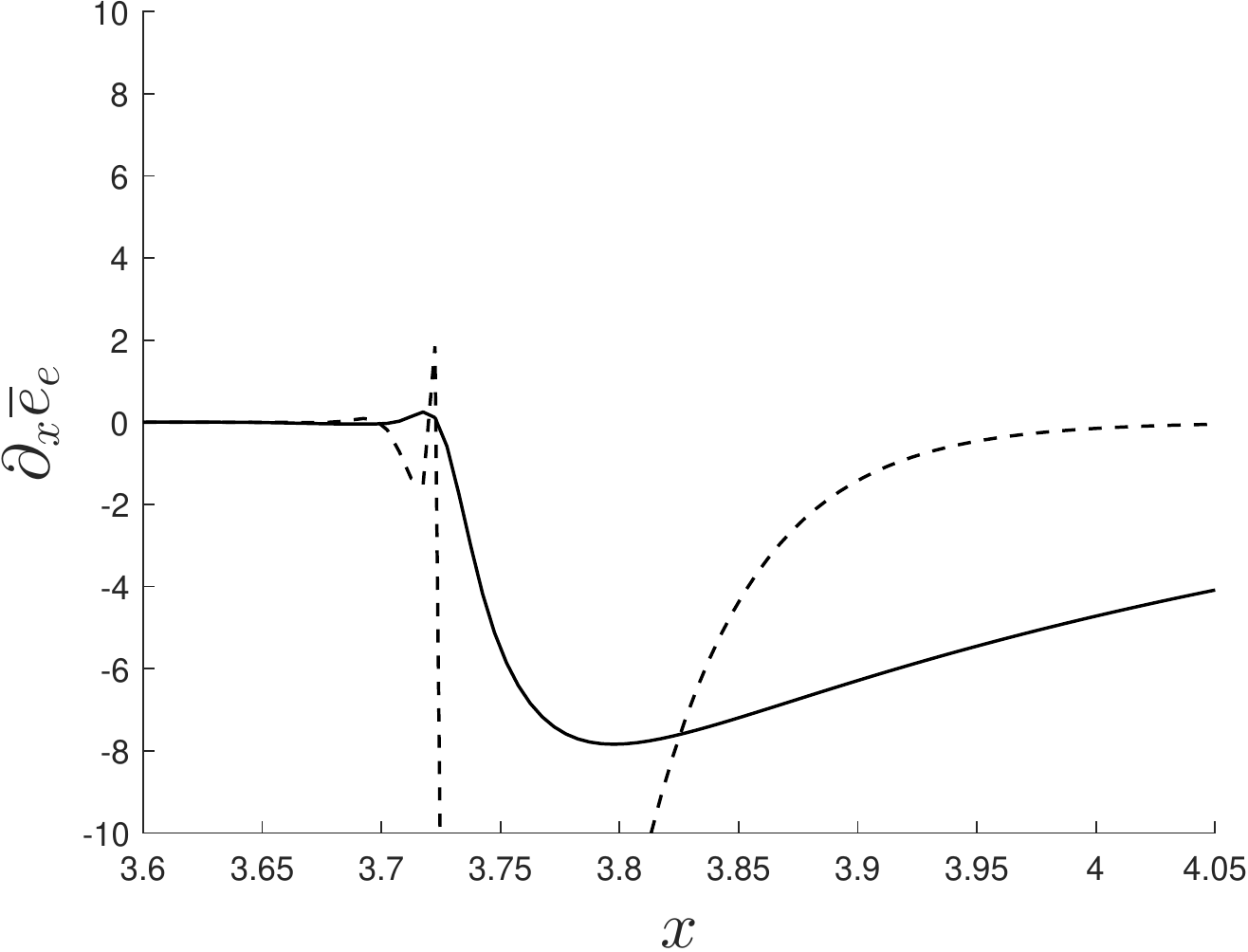}     
\centering   
\caption{Quantity $\dx \wenergiee$ for case \caseHD (full line) and case \caseWD (dashed line) at $t=t_f=1$.}
\label{fig:grad}
\end{figure}			

\Cref{fig:grad} represents the gradient of the electron energy $\dx\wenergiee=(\gamma {-}1)\dx\wtempe$ close to the discontinuity for the cases \caseHD and \caseWD at $t=t_f=n_f \Dt$, where $n_f$ is the total number of iterations at the final time $t_f$. For each cell $1\leq j\leq N$, the gradient is computed by means of a centered finite difference formula
\begin{equation*}
\gradwenergieejnf =\frac{\wenergieejpnf-\wenergieejmnf}{2\Dx}.
\end{equation*}
One can see that in the case \caseHD, the gradient is small whereas in the case \caseWD, a numerical artefact appears in the discontinuity. If the equality $\wtempe'(0^+)=\wtempe'(0^{-})$ is verified numerically, one gets the proper traveling wave with the right jump condition. If not, an artefact is produced in the discontinuity, due to the poor resolution of the gradient in the discontinuity, resulting in an artificial numerical shock.	
		
\begin{figure}[ht]
\centering
\includegraphics[scale=0.5]{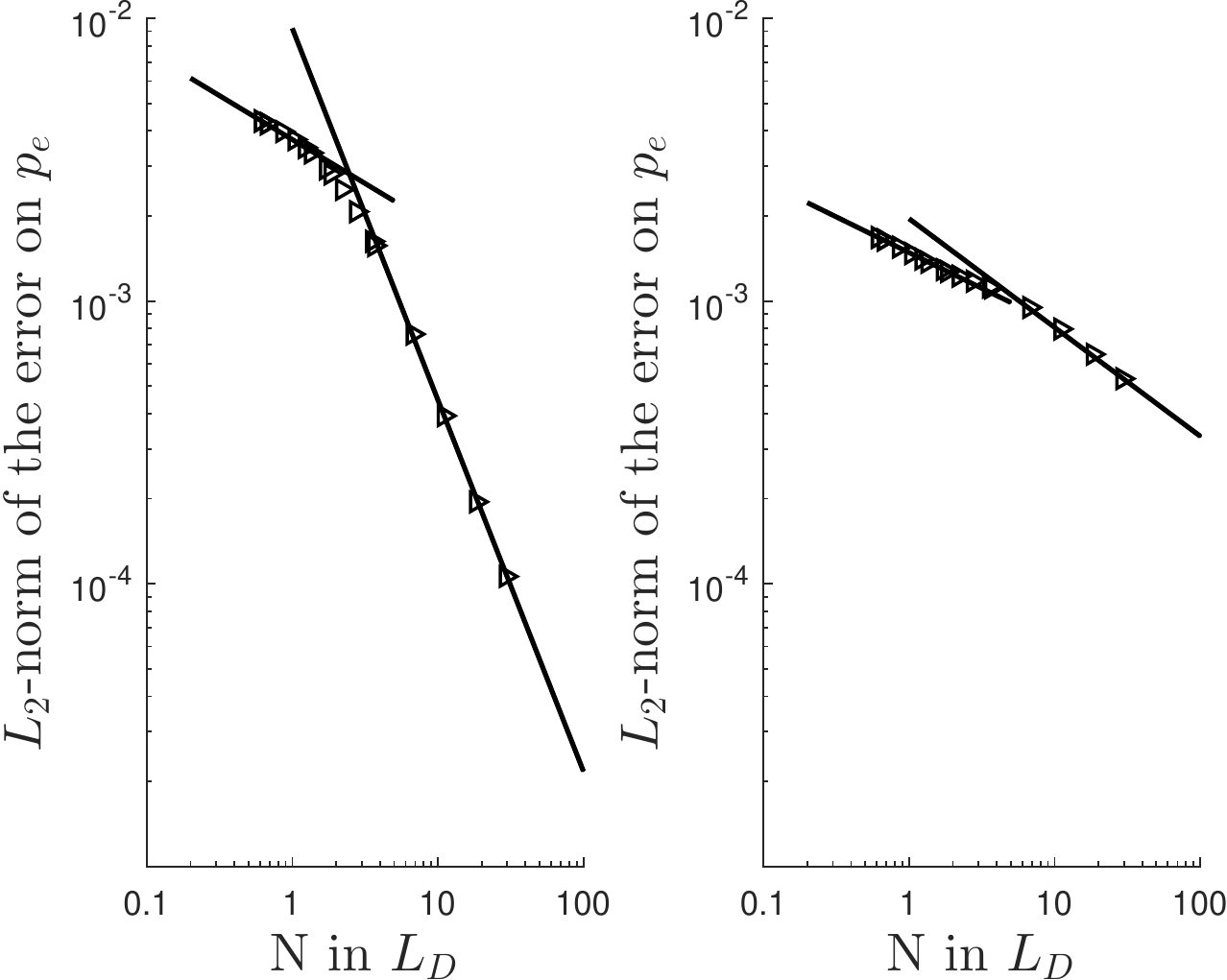}     
\includegraphics[scale=0.5]{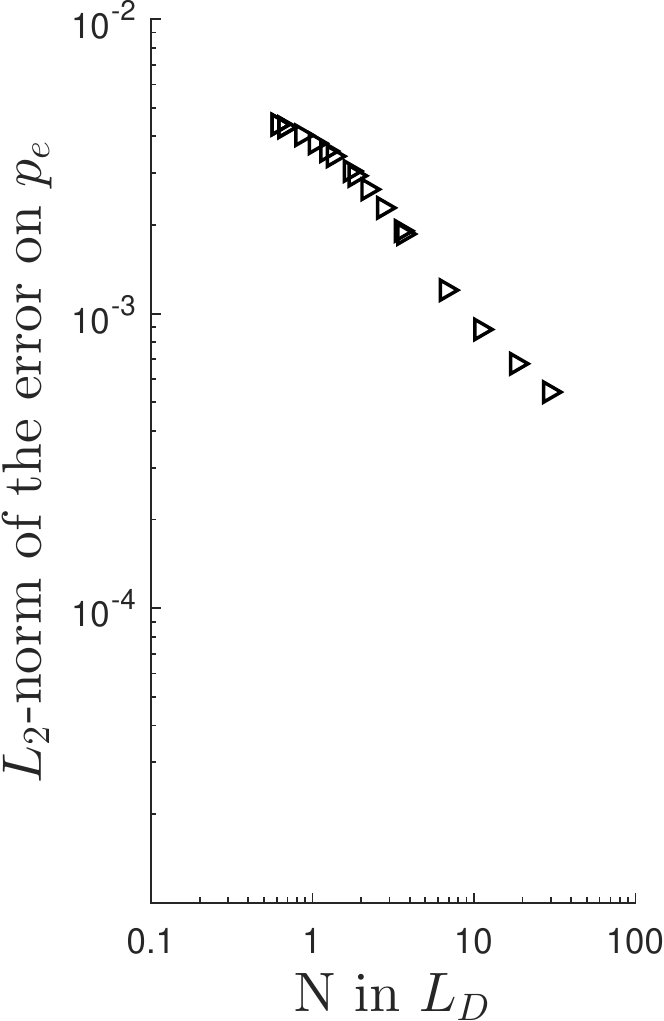}    
\caption{$L_2$-norm of the error on $\wpre$ with respect to the number of nodes in $\LD$. Left: downstream domain contribution (slopes of the lines: $0.3324$ and $1.314$); middle: upstream domain contribution (slopes of the lines: $0.2541$ and $0.3846$); right:  full domain contribution.}
\label{fig:L1norm}
\end{figure}

\Cref{fig:L1norm} shows the error in $L_2$-norm of $\wpre$ in function of the number of nodes in $\LD$, downstream and upstream of the shock. For the two areas studied, the dynamics of the $L_2$ norm of the error on $\wpre$ is the same: when the resolution of the wave is increasing, this norm is decreasing. Downstream of the shock, one can identify two dynamics for the $L_2$ norm. It can also be noted that the dynamic is changing when the error in $L_2$ norm   downstream of the shock becomes greater than upstream of the shock. Moreover, it is exactly at these resolutions that we begin to observe the appearance of an artificial shock.

The conclusive remark of these numerical experiments is the following. 
If the gradients in the discontinuity are well resolved (that is the case when the spatial mesh is fine enough), then,                                                                                                                                                    the traveling wave is well captured; 
if not, an artificial numerical shock is produced 
and, in this case, the numerical dissipation is responsible of the dynamic of the traveling wave. While the wave is regularized, having a non-linear model with diffusion coming from physics implies that the conditions for solving the wave are identified. On the contrary, in systems based on general nonconservative hyperbolic equations, it is difficult to clearly identify how numerical dissipation impacts the resolution of the wave \cite{abgrall}, \cite{driollet}.

Since the conditions for capturing the traveling wave have been identified, one wants to improve these results and build a new way for discretizing the nonconservative product allowing us to capture properly the traveling wave, even in weakly resolved cases. This is particularly relevant in solar physics. Indeed, it will be seen in \Cref{sec:sec5}  that the structure of the traveling wave in the sun chromosphere conditions corresponds to a case where the characteristic length $\LD$ is very small compared to the characteristic length $\LT$. Consequently, using the presented standard scheme with the standard discretization of the nonconservative product, one would need a lot of nodes in order to capture properly the traveling wave.
 
\subsection{Specific treatment of the nonconservative product}

In this section, we develop an original method for discretizing the nonconservative term $\Fjncdeux$. The idea is to express in a discretized sense the compatibility conditions  \cref{eq:discontinuity}  for the discontinuity, and deduce a new expression of the nonconservative term $\Fjncdeux$ in order to satisfy these conditions. In that way, we aim at capturing the traveling wave even when the gradients are not fully resolved. 

We consider the general scheme written in term of finite volumes \cref{eq:generalscheme} and 
specify the two vectorial coordinates by means of an index 
\begin{equation}
\left\lbrace
\begin{aligned}
&\Ujnpun-\Ujnun+\tfrac{\Dt}{\Dx}\bigl(\Fjcpun-\Fjcmun\bigr) = 0+\tfrac{\Dt}{\Dx}\bigl(\FjDpun-\FjDmun\bigr),\\
&\Ujnpdeux-\Ujndeux+\tfrac{\Dt}{\Dx}\bigl(\Fjcpdeux-\Fjcmdeux\bigr) = \Fjncdeux+\tfrac{\Dt}{\Dx}\bigl(\FjDpdeux-\FjDmdeux\bigr).
\end{aligned}	
\right.
\end{equation}

According to the relation found in \Cref{eq:discontinuity}, the jump of the gradient of $\wpre$ in the discontinuity is the same in the  electron mass and electron thermal  energy equations. In the discretized sense, one can simply link these two equations by multiplying the electron mass equation by a temperature \(\tempejn=(\gamma-1)\Ujndeux/\Ujnun\):
\begin{equation}
\left\lbrace
\begin{aligned}
&
\tempejn(\Ujnpun-\Ujnun)
+\tempejn\tfrac{\Dt}{\Dx}\bigl(\Fjcpun-\Fjcmun\bigr)
= \tempejn\tfrac{\Dt}{\Dx}\bigl(\FjDpun-\FjDmun\bigr),
\\
&\Ujnpdeux-\Ujndeux+\tfrac{\Dt}{\Dx}\bigl(\Fjcpdeux-\Fjcmdeux\bigr) = \Fjncdeux+\tfrac{\Dt}{\Dx}\bigl(\FjDpdeux-\FjDmdeux\bigr).
\end{aligned}
\right.
\end{equation}
According to \Cref{eq:discontinuitysimplified}, the derivative of the temperature $\wtempe$ is continuous in the discontinuity so  we have:\begin{equation}
\tempejpn-\tempejn = \tempejn -\tempejmn \iff 
\tempejpn-2\tempejn +\tempejmn = 0.
\label{eq:tempcon}
\end{equation}
Consequently, the second-order terms of the electron thermal energy equation can be simplified using \Cref{eq:diffusiveflux } as
\begin{equation}
\FjDpdeux-\FjDmdeux=\frac{D\gamma}{\Dx}\left(\Ujndeux[j+1]-2 \Ujndeux[j]+\Ujndeux[j-1] \right),
\end{equation}
and coupled to second-order terms of the  electron mass \Cref{eq:diffusiveflux },  defined as
\begin{equation}
\FjDpun-\FjDmun=\frac{2 D\left(\gamma-1 \right)}{\Dx}\left[\frac{\left(\Ujndeux[j+1]-\Ujndeux[j]\right)}{\left(\tempejpn+\tempejn\right)} - \frac{\left(\Ujndeux[j]-\Ujndeux[j-1]\right)}{\left(\tempejn+\tempejmn\right)} \right],
\label{eq:form}
\end{equation}
Besides, the time derivative terms are not playing a role in the compatibility  \Cref{eq:discontinuity}. We couple the  electron mass  equation with the  electron thermal energy equation,  providing a new expression for the nonconservative product. Multiplying the first line of \Cref{eq:tempcon} by $\gamma/(\gamma-1)$, substracting by the second line of \Cref{eq:tempcon}, and neglecting the temporal derivative, one obtains

\begin{multline}
\Fjncdeux = \tfrac{\Dt}{\Dx}\bigl(\Fjcpdeux-\Fjcmdeux\bigr)-\tfrac{\gamma}{\gamma-1}\tempejn \tfrac{\Dt}{\Dx}\bigl(\Fjcpun-\Fjcmun\bigr)\\
-\tfrac{\Dt}{\Dx}\bigl(\FjncDp-\FjncDm \bigr),
\label{eq:formulation}
\end{multline}
where the second-order terms of \Cref{eq:formulation} are given by
\begin{equation}
\FjncDpm=\tfrac{1}{2}\tfrac{\gamma}{\gamma-1}\bigl(\tempejpmn-\tempejn\bigr)\FjDpmun.
\end{equation}

This expression of the nonconservative product \smash{$\Fjncdeux$} verifies \Cref{eq:discontinuity} in the discretized sense, that is to say: 1- the continuity of $\wtempe$ in the discontinuity of the traveling wave verifying \Cref{eq:tempcon}, 2- the conditions on the jump of the gradient of $\wpre$ in the discontinuity. 

However, the expression found in \Cref{eq:formulation} for $\Fjncdeux$ makes the global scheme not consistent with system 
\eqref{eq:couple2}, even if it is able to resolve the traveling wave in the weakly discretized case. It is thus necessary to 
add correction terms to recover consistency in the well-resolved case.

Thus, we want to build a numerical scheme in order to get proper compensations of the different terms in discontinuities in order to verify the compatibility condition \Cref{eq:discontinuity} in the discretized sense, and at the same time, to add correction terms and recover consistency in order to capture properly the regular parts of the traveling wave. To this aim, we focus only on the first order terms of the specific expression of the nonconservative terms in \Cref{eq:formulation}. Then, we add first order correction terms in order to be consistent with the original system \eqref{eq:couple2}. Several alternatives are possible to define the correction terms. In this work, we have mainly focused on the correction terms, which involve the gradients of $\wrhoe\wenergiee$ and $\wenergiee$ (or $\wtempe$). By adding first order correction terms, the expression of the nonconservative product \Cref{eq:formulation} reads:
\begin{equation}
\Fjncdeuxcorrection = \Fjncdeux 
-\Dt~\wvitessejn\frac{\Ujndeux[j]-\Ujndeux[j-1]}{\Dx}
+\Dt\frac{\gamma}{\gamma-1} \wvitessejn  \tempejn \frac{ \Ujnun[j] - \Ujnun[j-1] }{\Dx}.
\label{eq:formulationfinal}
\end{equation}

\begin{figure}
\centering
\includegraphics[scale=0.5]{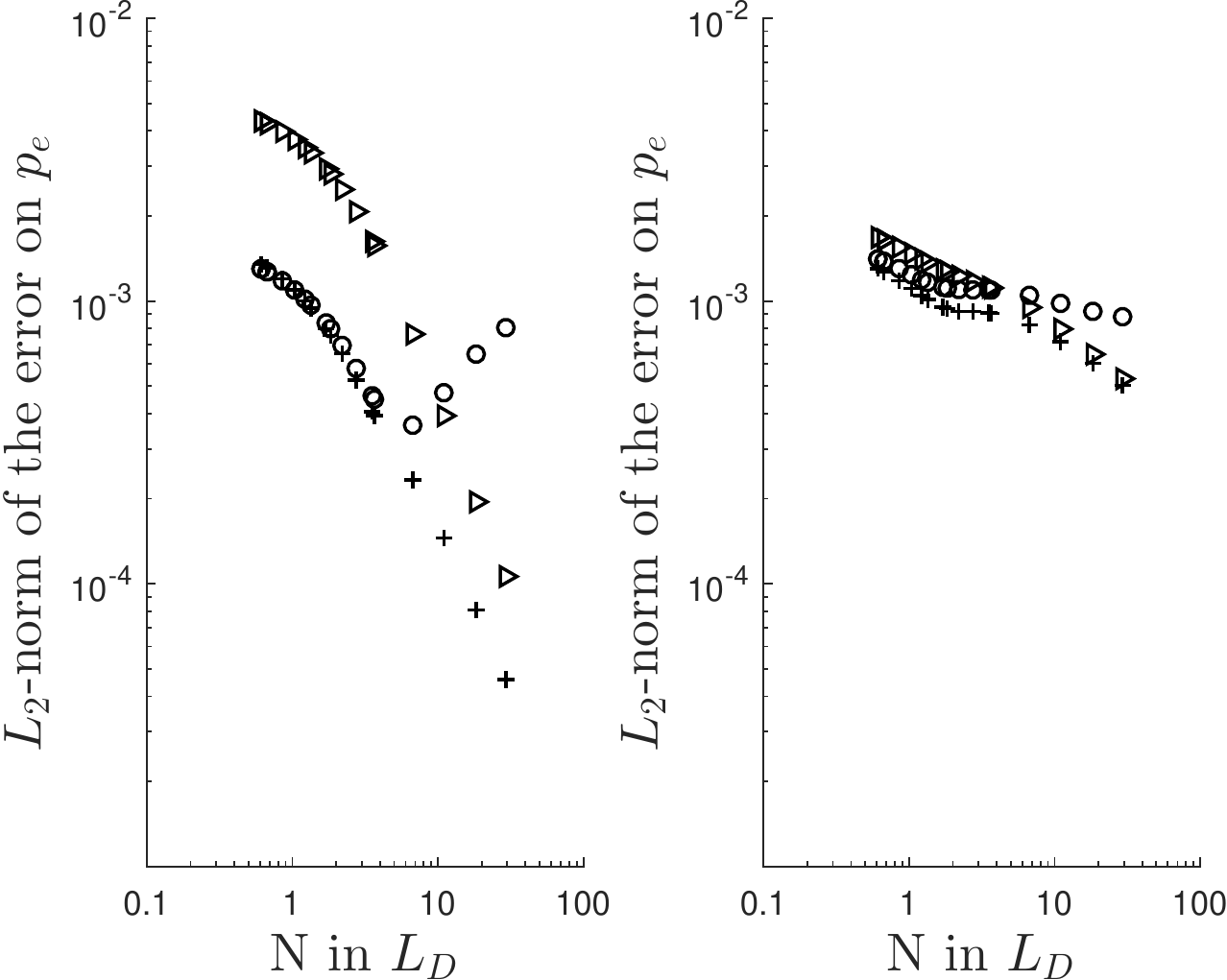}     
\includegraphics[scale=0.5]{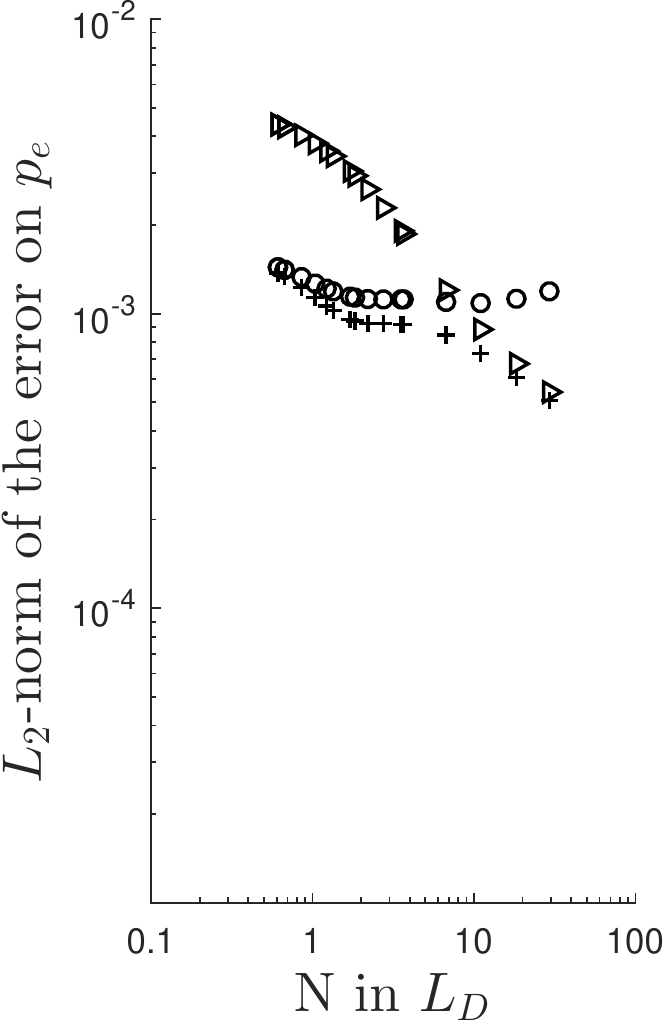}    
\caption{$L^2$-norm of the error on $\wpre$ with respect to the number of nodes in $\LD$, $\triangleright$~standard discretization, $\circ$~discretization without correction terms from \Cref{eq:formulation}, and $+$~discretization with correction terms from \Cref{eq:formulationfinal}. Left: downstream domain contribution; middle: upstream domain contribution; right:  full domain contribution.}
	\label{fig:norm2}
\end{figure}

It is important to make sure that these additional correction terms are not playing any role in discontinuities, which could violate \Cref{eq:discontinuity} in the discretized sense. We limit the impact of these additional terms in the discontinuities and one can show that the dynamics of the wave is not depending on the choice of the cut-off.

\Cref{fig:norm2} represent the $L^2$ norm of $\wpre$ in function of the number of nodes in $\LD$ for three different way of discretizing the nonconservative product $\Fjncdeux$. Results are presented for three methods: the standard way of discretizing $\Fjncdeux$ described in \Cref{sec:sec3}, without correction terms of \Cref{eq:formulation}, and with correction terms $\Fjncdeuxcorrection$ of \Cref{eq:formulationfinal}. The formulation of the nonconservative product $\Fjncdeux$ without correction terms is working well in the weakly discretized case for the traveling wave test cases. However, in the regularized case, an additional numerical shock is appearing and the $L_2$ norm is increasing with the number of nodes in $\LD$, because the scheme is not consistent. By adding correction terms and using the formulation of the nonconservative term $\Fjncdeuxcorrection$ from \Cref{eq:formulationfinal}, we have built a scheme which is able to capture the traveling wave in both the highly- and coarsely-resolved cases.

We could have built this numerical scheme because of the thorough understanding of the traveling wave linked to a good structure of the diffusion, as well as system allowing us  to derive compatibility equations in the discontinuity.

\section{Application to solar physics}
\label{sec:sec5}

In this section, we apply the previous development to a test case chosen so as to reproduce typical scales from sun chromosphere conditions. We study the ability of our scheme to resolve shock solutions in such conditions and design a specific numerical strategy based on the new scheme in order to cope with the nonconservative term. A 3-wave is considered by using the system \cref{eq:couple2} with non-dimensional quantities for building the traveling wave in the sun chromosphere conditions.

First, we use atmospheric parameters from the model C of Vernazza \textit{et al.} \cite{vernazza} where the values of these parameters are given at $52$ depth in the atmosphere from the low corona to the photosphere. For the purpose of the work, we have focused on the photospheric level at the heigth  $h=0\, km$. We consider $\wpreR=\wpri[\heavy][\droit]$ at the right state, with only two species: electrons and protons $H^+$ as heavy particles. The transport coefficients $D$ and $\lambda$ are computed using third-order Sonine polynomials approximation based on a spectral Galerkin method used in \cite{magin04,wargnier2018} considering local thermodynamic equilibrium for the fully ionized gas. 

Then, we non-dimensionalize these quantities with reference quantities. The characteristic length of diffusion $L_0=\LD$ is chosen as the reference length. The density of heavy particle is chosen as the reference density $\rho_0$. The reference velocity $v_0$ is the Alfv\'en velocity defined as
\(v_0=B_0/\sqrt{\mu_0 \rho_0}\)
where $B_0$ is the reference magnetic field, chosen as $B_0 = 100$ G, and $\mu_0$ the vacuum permeability.

\begin{table}[ht]
\centering\footnotesize
  \caption{Reference quantities at the photospheric level}
  \label{tab:tablesolar3}
\grandtraittop
\renewcommand*{\arraystretch}{1.1}
\begin{tabular}{@{}p{2cm}p{2cm}p{2cm}p{1.2cm}p{1.2cm}p{2cm}@{}}
$\rho_0$ ($kg.m^{-3}$) & $L_0$ ($m$)& $v_0$ ($m.s^{-1}$)& $T_0$ ($K$)& $P_0$ ($Pa$)&$n_0$ ($m^{-3}$)\\ \hline
 $1.873\times 10^{-4}$ & $1.747\times 10^{-6}$&$6.518\times 10^{2}$&6420&9927.42&$1.12\times 10^{23}$
\end{tabular}
\grandtraitbottom
\end{table}

Finally, after non-dimensionalizing the governing equations with reference quantities from \Cref{tab:tablesolar3}, the traveling wave is investigated using values from \Cref{tab:tablesolar} and \Cref{tab:tablesolar2}. We have chosen a number of nodes of either $N=1000$ or $N=5000$ and a length of the domain $L/L_0=2{\times}10^{5}$.

\begin{table}[ht]
\centering\footnotesize
  \caption{Right and left states at infinity of the traveling wave}
  \label{tab:tablesolar}
\grandtraittop
\renewcommand*{\arraystretch}{1.1}
\begin{tabular}{@{}p{2cm}p{1.2cm}p{2cm}p{1.2cm}p{1.2cm}p{1.2cm}@{}}
&$\wrhoH$ & $\wrhoe$ & $\wpression$&$\wpre$&$\wvitesse$\\ \hline
right state $R$ & 1 &$5.44{\times}10^{-4}$ &0.5974&0.2987&0.07\\
left state $L$ &1.6962&$9.23{\times}10^{-4}$&1.5&0.9454&0.6787
\end{tabular}
\grandtraitbottom
\end{table}

\begin{table}[ht]
\centering\footnotesize
  \caption{Diffusion coefficients and related typical lengths}
  \label{tab:tablesolar2}
\grandtraittop
\renewcommand*{\arraystretch}{1.1}
\begin{tabular}{@{}p{1.2cm}p{2cm}p{1.2cm}p{1.2cm}p{1.2cm}@{}}
D & $\kappa^\droit$ & $\LD/L_0$&$\LT/L_0$&$L/L_0$\\ \hline
 10.7853 &121 970.96&1&11309 &200 000
\end{tabular}
\grandtraitbottom
\end{table}

After initializing the traveling wave, three numerical schemes are compared. The first scheme, denoted \textbf{scheme A}, is the standard scheme based on a standard discretization of the nonconserva\-tive product introduced in \Cref{sec:sec4}, where the timestep, denoted  $\Dt_F$, is limited by a  Fourier stability condition, thus involving the largest diffusion coefficient, that is the electron thermal diffusivity $\kappa^\droit$. The second scheme, denoted \textbf{scheme B}, is based on the formulation of the nonconserva\-tive product described in \Cref{eq:formulationfinal}, where the timestep is also limited by the same Fourier stability condition. The third scheme, denoted \textbf{scheme C}, is based on the formulation of the nonconserva\-tive product defined in \Cref{eq:formulationfinal}, using an operator splitting approach based on a second-order Strang formalism in order to separate the convection and diffusion operators \cite{strang,duarte,Duarte14,theseduarte}. The idea is to not be limited by the small timesteps $\Dt_F$ imposed by the Fourier stability condition, due to the electron thermal diffusivity. Concerning the \textbf{scheme C}, there are two possibilities for the operators : one can 1- gather diffusive terms and the nonconservative product, or 2- gather convective terms and the nonconservative product. In this work, we have focused on the second case. Indeed, according to \Cref{eq:formulationfinal}, the expression of the nonconservative product depends on the thermal energy and density convective fluxes, which makes the second case a rather more natural choice. Besides, the computational time is drastically shorter in this case, since the nonconservative product is integrated only one time during the convective timestep, whereas in the other case, the nonconservative product is integrated several times  during the dissipative timestep.

The operators are splitted: one operator $X$ corresponds to convective terms and the nonconservative product defined by \Cref{eq:formulationfinal}, where the convective timestep, called $\Dt$, is simply limited by a CFL condition; an other operator $Y$ regroups diffusion terms, where the timestep $\Dt_F$ is computed based on the Fourier condition and integrated over several sub-timesteps in order to reach the convective timestep. The general scheme is summarized as follows:
\begin{equation*}
U^{n+1}=Y^{\frac{\Dt}{2}} X^{\Dt}Y^{\frac{\Dt}{2}} U^{n}.
\end{equation*}
In the proposed schemes, the values used for the convective timestep $\Dt$ and the diffusive timestep $\Dt_F$ are presented  in \Cref{table:timestep}. In order to perform the comparison between the timesteps used, we have  compared them to the convective timestep $\Dt=C\times\Dx/\max(\wvitesse^\droit+\cR,\wvitesse^\gauche+\cL)$ for $N=1000$ and $N=5 000$, where the Courant number is $C=0.2$.

\begin{table}[ht]
	\centering\footnotesize
	\caption{Timesteps used for the three schemes for $N=1000$ and $N=5000$}
	\label{table:timestep}
	\grandtraittop
	\renewcommand*{\arraystretch}{1.1}
	\begin{tabular}{@{}p{1.2cm}p{2cm}p{2cm}@{}}
		N& $\Dt$ & $\Dt_F$\\ \hline
		$1000$&$2.233\times 10^{1}$& $4.095\times10^{-1}$ \\
		$5000$&$4.466$ & $1.64\times10^{-2}$ \\
	\end{tabular}
	\grandtraitbottom
\end{table}

\begin{figure}
\centering
\subfloat[$N=1000$]{\label{fig:casesolarsplit}\includegraphics[scale=0.25]{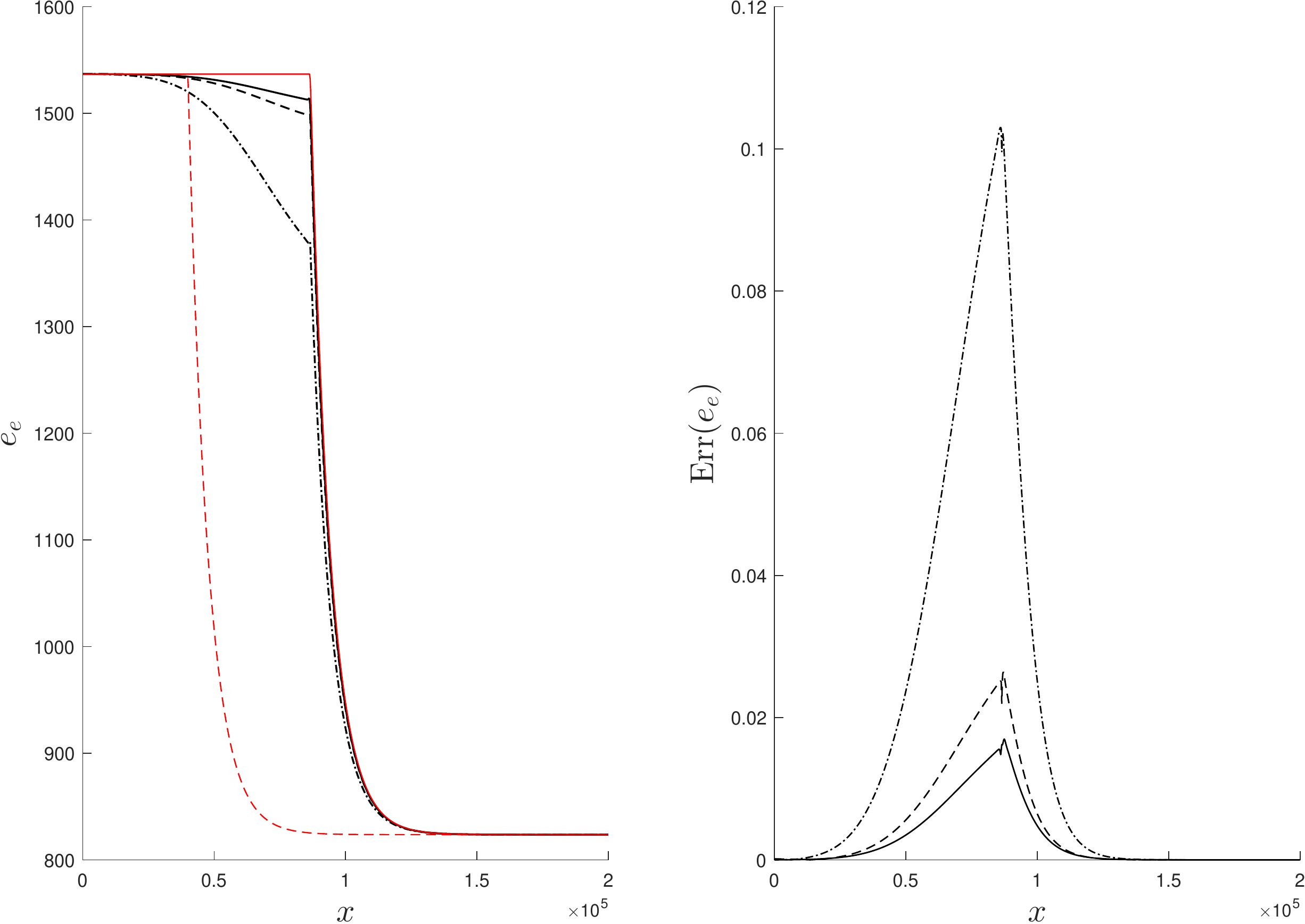}}
\hfill
\subfloat[$N=5000$]{\label{fig:casesolarsplit2}\includegraphics[scale=0.25]{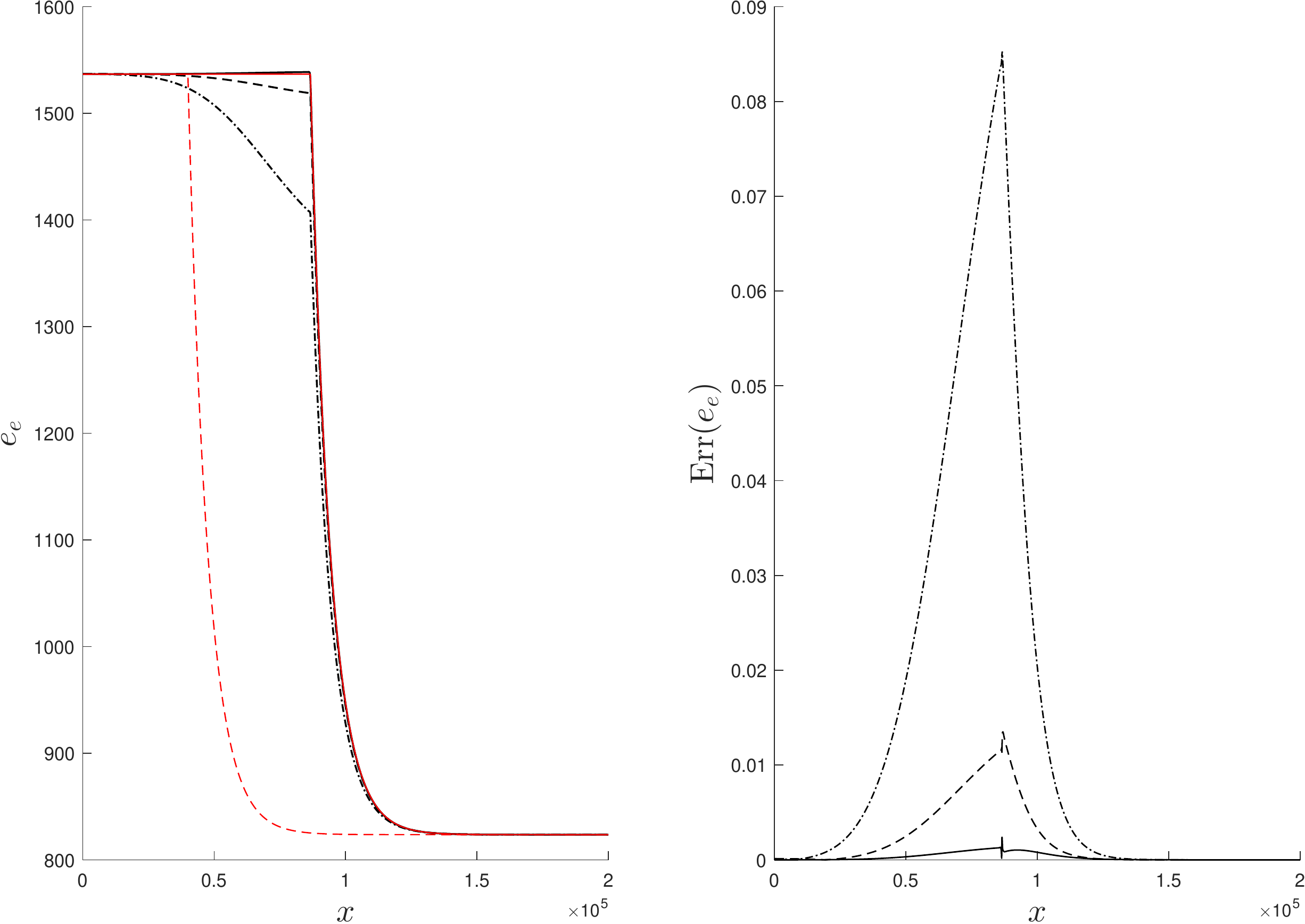}}
\caption{Electron energy ($\wenergiee$) and relative error (Err) for the solar test case based on the values from \Cref{tab:tablesolar} and \Cref{tab:tablesolar2}. Exact solution at $t=0$ (red dashed line) and $t=30 000$ (red full line). Numerical solution  for \textbf{scheme A} (semi-dashed line),  \textbf{scheme B} (dashed line), and \textbf{scheme C} (full line), at the final time $t=30 000$.}
\end{figure}

Results are presented in \Cref{fig:casesolarsplit} and \Cref{fig:casesolarsplit2}
for the electron energy $\wenergiee$, comparing the three schemes, at the final time $t=t_f=30 000$, for $N=1000$ and $N=5000$. 
In the sun chromosphere conditions, the characteristic scales are such that $\LT\gg\LD$, because the electron diffusivity is much higher than the electron diffusion coefficient in such conditions. In the test case, the smallest spatial scale to be resolved is the length $\LD$, which is the characteristic scale related to the resolution of the traveling wave. By fixing the number of nodes $N$ and the length of the domain $L$ based on \Cref{tab:tablesolar2}, the test case can be identified as a very weakly-resolved traveling wave test case. 

In \Cref{fig:casesolarsplit} and \Cref{fig:casesolarsplit2}, the standard scheme exhibits artificial numerical shock since the smallest scale is not properly resolved and small timesteps $\Dt_F$ have to be used, as expected. Switching to a proper treatment of the nonconservative term allows to reduce by a factor of $9$, for $N=5000$, the amplitude of the error on the electron temperature and to reduce drastically the artificial numerical shock, thus leading to a satisfactory level of resolution.
However, using the new scheme based on Strang splitting operator techniques combined to the new formulation of the nonconservative product (\textbf{scheme C})  leads to an additional improvement of the resolution of all the scales of the traveling wave. The traveling wave can be well captured, while using splitting timesteps of the order of the convective CFL stability limitation, thus leading to a minimal amount of numerical dissipation in the convective step. 
In fact, based on the results obtained in the previous section, a good approximation of the traveling wave obtained with the \textbf{scheme A} would require several nodes in the characteristic length $\LD$ thus leading to about a million nodes. In this context, the corresponding convective timestep and the diffusive timestep would be respectively $\Dt = 8.94\times 10^{-2}$ and $\Dt_F=4.01\times 10^{-7}$ and the original scheme would become useless.

\begin{figure}
\centering
\includegraphics[scale=0.3]{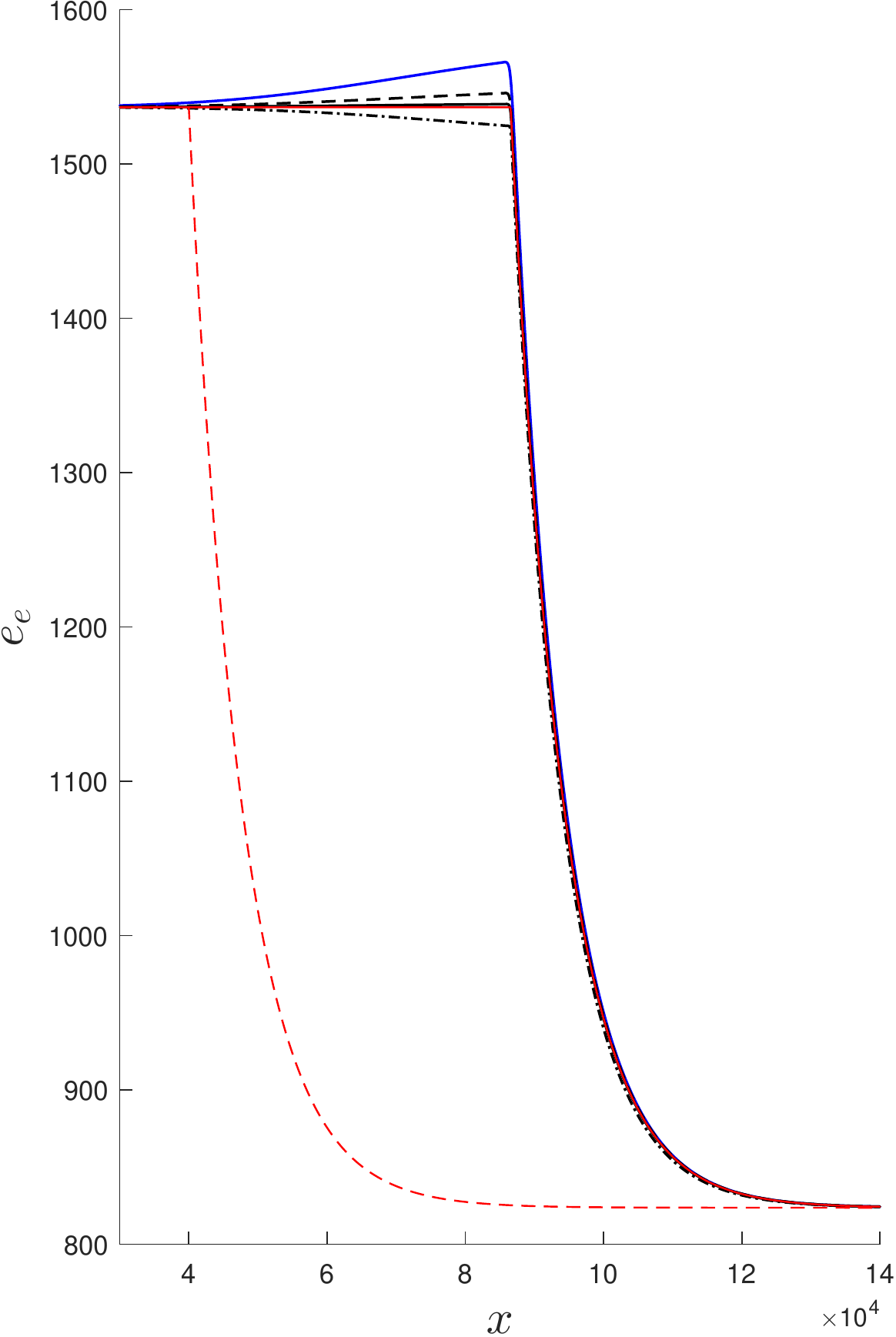}     
\caption{Electron energy ($\wenergiee$) and relative error (Err) for the solar test case based on the values from \Cref{tab:tablesolar} and \Cref{tab:tablesolar2}. Exact solution at $t=0$ (red dashed line) and $t=30 000$ (red full line). Numerical solution for \textbf{scheme C}  at the final time $t=30 000$, for $N=5000$ in the domain $\left[30000, 140 000\right]$ for several Courant number $C=5\times 10^{-2}$ (semi-dashed line),
$C=0.2$ (full black line),
$C=0.3$ (dashed line),
$C=0.4$ (full blue line).}
\label{fig:casesolarsplit3}
\end{figure}

In order to perform the analysis of the error generated by the splitting operation, several splitting timestep $\Dt$ have been tested for the presented test case, used in the \textbf{scheme C}. Results for several Courant number $C\in\{0.05,0.2,0.3,0.4\}$ are shown in \Cref{fig:casesolarsplit3}, for $N=5000$, at $t=30 000$. These results show that when the splitting timestep becomes too important (for Courant numbers $C= 0.3$ or $C= 0.4$), the traveling wave  is no longer captured with a high level of accuracy, and an additional numerical artefact is obtained. However, for a splitting timestep where the Courant number is $C= 0.2$, we notice that the numerical solution is shown to be optimal.

\section{Traveling wave for the fully coupled problem \cref{eq:fullprob}}
\label{sec:sec6}

In this paper, we have focused so far on the decoupled problem \cref{eq:couple1} and \cref{eq:couple2}. In the case of the fully coupled problem \cref{eq:fullprob}, the problem is in fact very similar to the decoupled problem and one can solve for traveling wave solutions as well. However, we get a numerical solution instead of a complete analytical solution. In this study, we also consider a 3-wave.

\begin{figure}[htbp]
\center{
\tikzset{
    fleche/.style = {->, thick},
    courbe/.style = {mark=none, very thick, samples=100},
    pente/.style = {mark=none, thick, dashed},
}
\begin{tikzpicture}[scale = .6]
\begin{axis}[
  grid=major,
  xtick={0, 1},
  xticklabels={0, $\LD^{*}$},
  ytick={0, .3, 1},
  yticklabels={0, ${\pre[\droit]}$, ${\pre[\gauche]}$},
  xlabel={$\xi=x-\sigma t$},
  ylabel={$\pre$},
  xlabel style={below},
  ylabel style={right},
  xmin=-5,
  xmax=5.5,
  ymin=-0.1,
  ymax=1.25,
  ]
\addplot+[courbe, draw=black, domain=-5:0] {%
	1%
	};
\addplot+[courbe, draw=black, domain=0:5] {%
	.3+.7*exp(-x)
	};
\addplot+[pente, draw=black, domain=0:1] {%
	1-x
	};
\end{axis}
\end{tikzpicture}
\hfill
\begin{tikzpicture}[scale = .6]
\begin{axis}[
  grid=major,
  xtick={0},
 xticklabels={0},
  ytick={.2,0.5, 1},
  yticklabels={${\energie^\droit,\vitesse^\droit,\rhoH[\droit]}$,${\energie^0,\vitesse^0,\rhoH[0]}$, ${\energie^\gauche,\vitesse^\gauche,\rhoH[\gauche]}$},
  xlabel={$\xi=x-\sigma t$},
  ylabel={${\energie,\vitesse,\rhoH}$},
  xlabel style={below},
  ylabel style={right},
  xmin=-5,
  xmax=5.5,
  ymin=-0.1,
  ymax=1.25,
  ]
\addplot+[courbe, draw=black, domain=-5:0] {%
	1%
	};
\addplot+[courbe, draw=black, domain=0:5] {%
	.2+.3*exp(-x)
	};
\end{axis}
\end{tikzpicture}
}
\caption{Structure of the traveling wave for the fully coupled problem \cref{eq:fullprob}}
\label{fig:schemafull}
\end{figure}
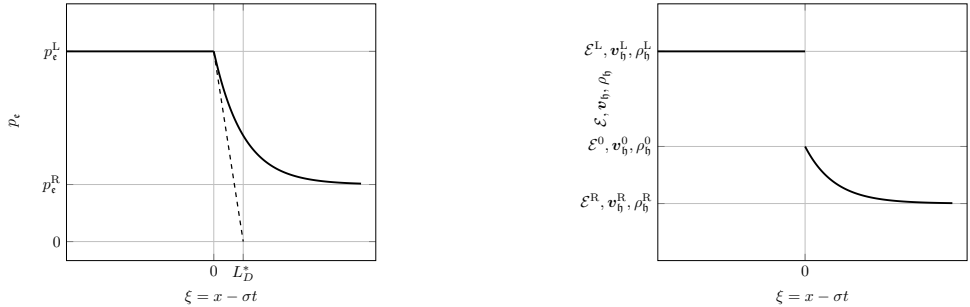

In the fully coupled system \cref{eq:fullprob}, we solve for a traveling wave where the structure corresponds to 1- a constant state $\gauche$, a weak discontinuity (smooth function with jump of derivative) and a regularization up to a constant state $\droit$ for the electron variables $({\pre,\rhoe})$ and 2- a constant state $\gauche$, a discontinuity connecting state $\gauche$ to an intermediary state $0$ and a regularization from state $0$  to a constant $R$ for the heavy variables (${\energie,\vitesse,\rhoH}$). The structure of the wave is represented in \cref{fig:schemafull}. The jump conditions are the Rankine-Hugoniot conditions, the thermal energy of electrons requires a numerical integration and the velocity jump is coupled to the weak jump of the electron variables. Actually, the structure of the wave for the heavy particles is very similar to the one identified by Zel\textsc{\char13}dovich and Raizer in \cite{zeldovich} in \S 3, where the role of the heat conduction on the structure of the shock wave in gases has been studied. 

In order to get an accurate estimation of the missing jump condition for the electron thermal energy of as a function of the Mach number $\MR$ as well as the structure of the traveling wave, we integrate the ordinary differential equations of \cref{eq:fullprob} using a Dormand-Prince (RKDP) method or DOPRI853 method \cite{DORMAND}. In the case of a 3-wave, we initialize the traveling wave from state $\gauche$ and, by numerical integration, compute the corresponding state $\droit$. In order to get an accurate estimation of the missing jump condition, a shooting method is used. The steps of the shooting method are the following: 1- we start the numerical integration of the traveling wave using the state $\gauche$ of the \textbf{case A} in  \cref{sec:num} of the decoupled problem, then 2- a state $\droit$ associated to the initial state $\gauche$ is computed, finally 3- a dichotomy is used by initializing different state $\gauche$ until a good approximation of the expected state $\droit$ is found. The missing jump condition of the fully coupled problem as a function of the Mach number can thus be obtained. 

\begin{figure}
\includegraphics[scale=0.4]{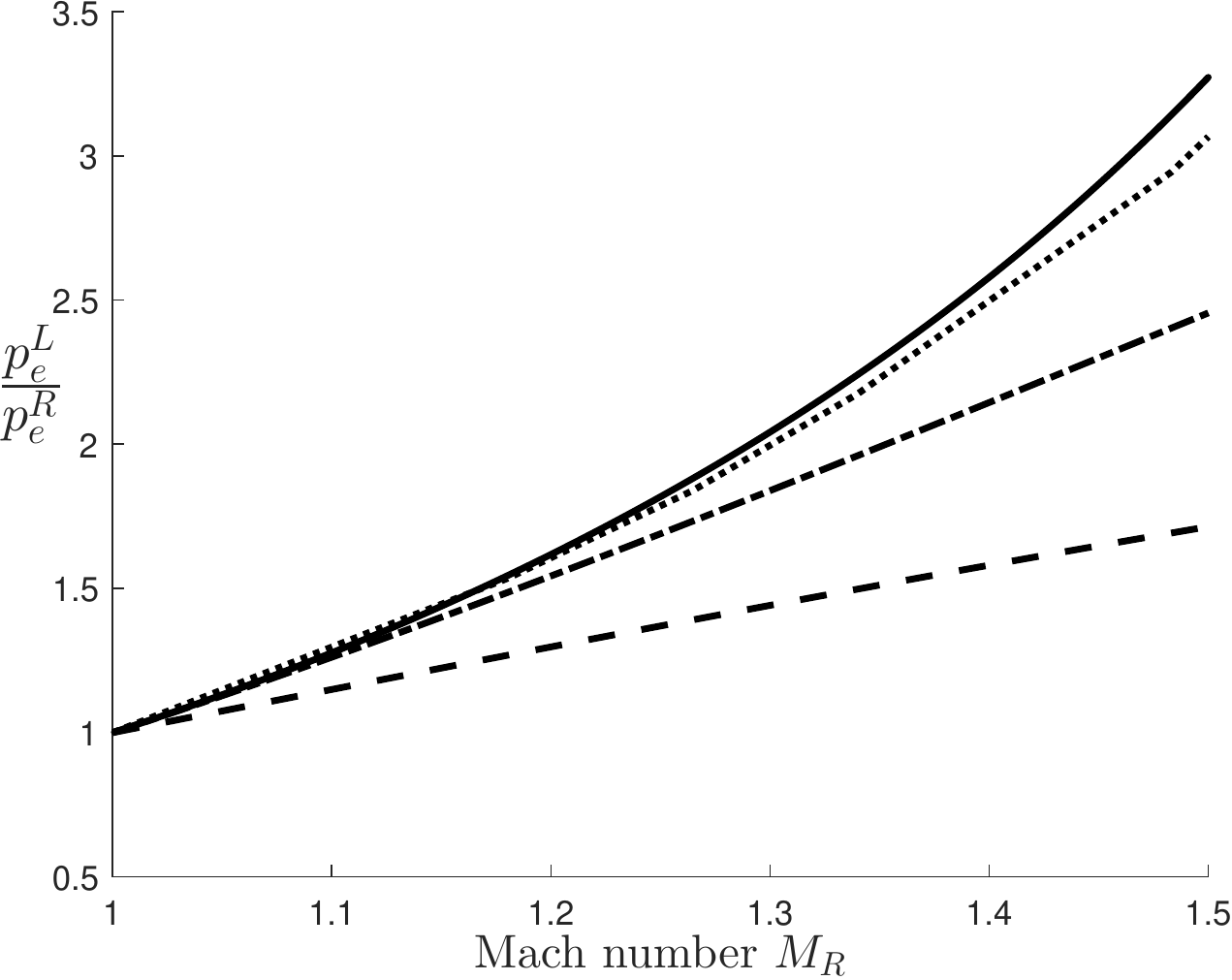} \hfill
\includegraphics[scale=0.4]{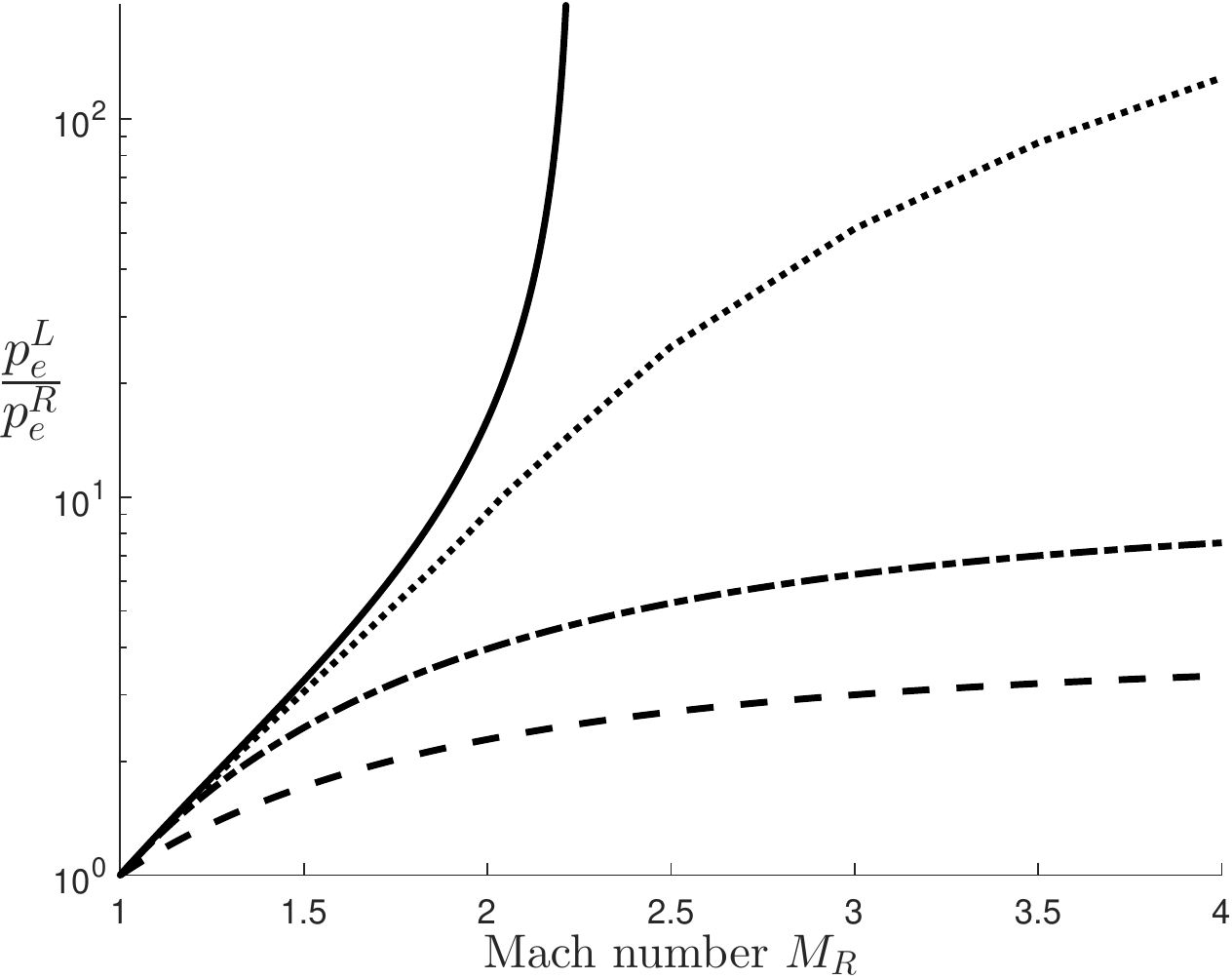} \hfill
\includegraphics[scale=0.4]{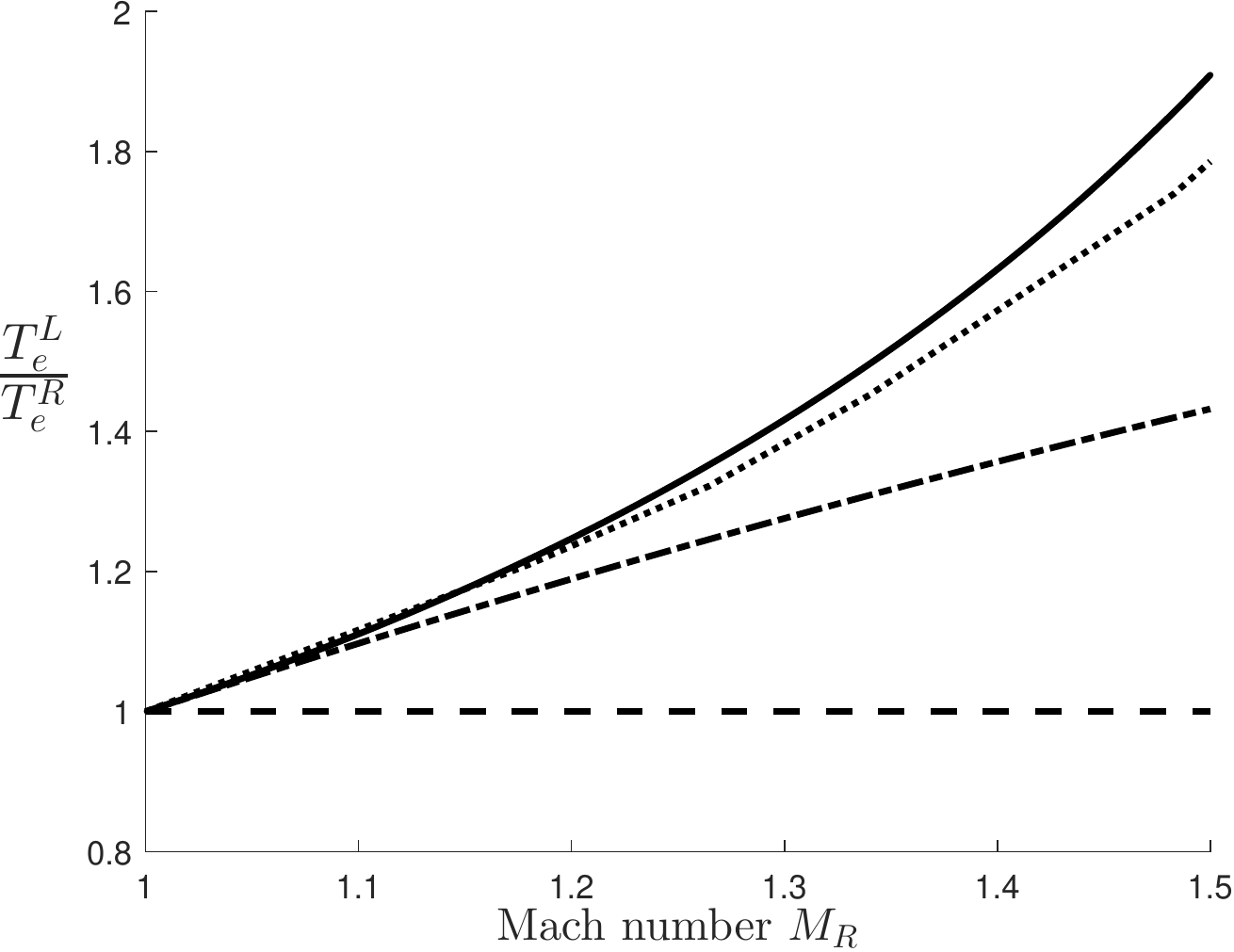} \hfill
\includegraphics[scale=0.4]{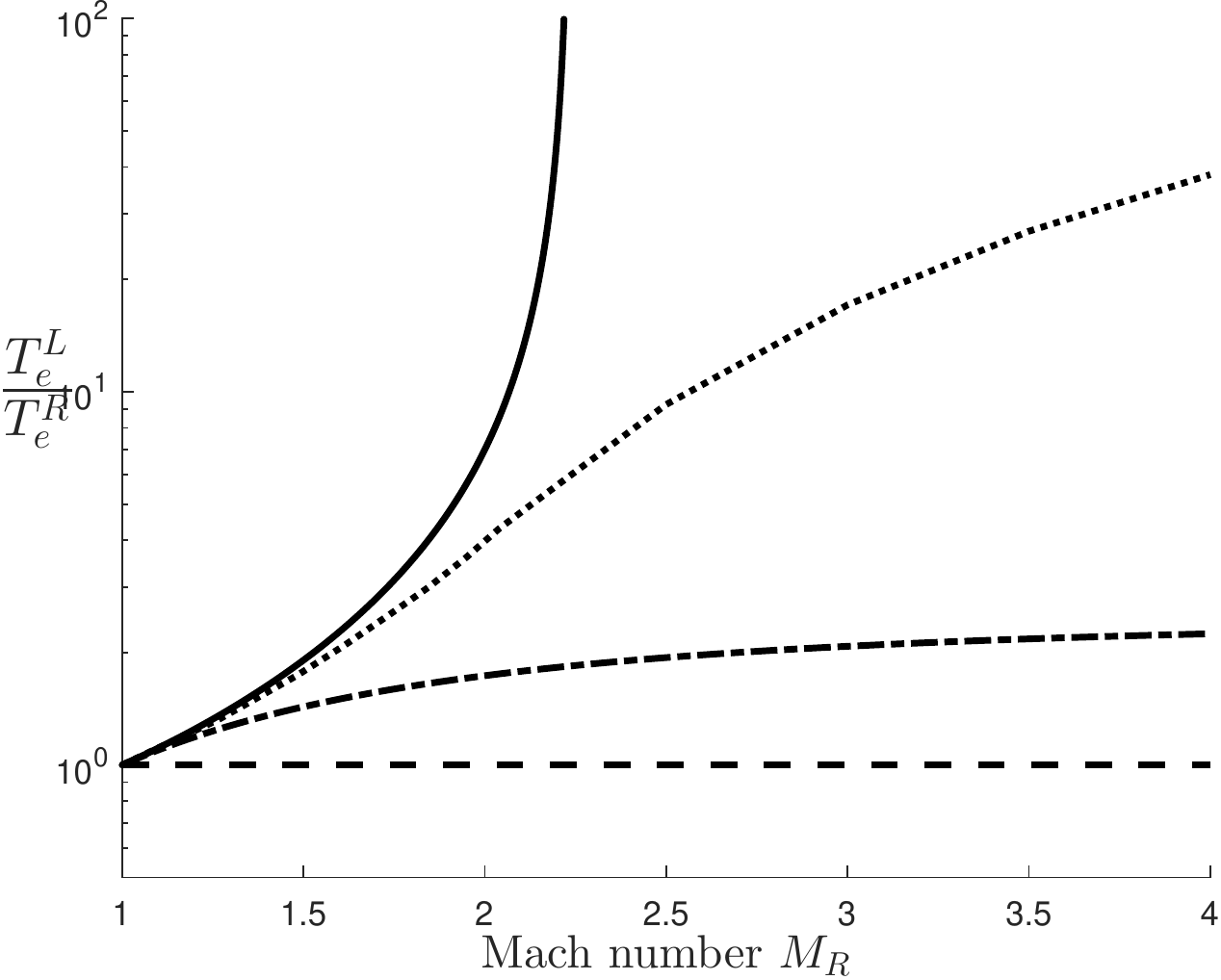}

\caption{Jump of $\pre$ and $\tempe$ as a function of the Mach number $\MR$.
\textbf{In full line}: the jump from the decoupled system \cref{eq:couple2}, in \textbf{ semi dashed-line}: model $\mathcal{M}^{\text{ent}}$ from the system \cref{eq:entropy}, in \textbf{dashed line}: model $\mathcal{M}^{\text{src}}$ from the system \cref{eq:source}, in \textbf{dotted line}: jump from the fully coupled system from \cref{eq:fullprob}} 
\label{fig:hmj}
\end{figure}

The results of the numerical integration are presented in \Cref{fig:hmj}. The estimated jump condition for the fully coupled problem is compared with the jump conditions from the decoupled problem \Cref{eq:jumppe}, \Cref{eq:jumppee} from the model $\mathcal{M}^{\text{ent}}$, and \Cref{eq:jumppes} from the model $\mathcal{M}^{\text{src}}$. The jump conditions of the decoupled problem \Cref{eq:jumppe} give a very good approximation of the jump conditions of the fully coupled problem in a reasonable Mach number range $\MR$ close to one. Besides, in the fully coupled case, no singularities have been observed for the jump condition of $\pre$ and $\tempe$ for the whole range of Mach number. The results show that the jump conditions from the conservative models $\mathcal{M}^{\text{ent}}$ and $\mathcal{M}^{\text{src}}$ clearly underestimated the post-shock temperature.

Finally, if having an analytical expression of the traveling wave for the fully coupled problem is not possible, relying on the same strategy designed in the study of the decoupled problem, we are able to analyze the fully coupled case. By integrating the ordinary differential equations of the fully coupled problem \cref{eq:fullprob}, we get the structure of the traveling wave as well as an evaluation of the missing jump condition of the internal energy of electron.  

From this study, two conclusions can be drawn: 1- Relying on the missing jump condition proposed in the literature through various approximation yields a very poor approximation of the 
effective jump conditions, even in a Mach number range close to one and the present study allows to derive the physically sound jump conditions, 2- Focusing on the decoupled problem, at least in a reasonable Mach number range around one, is fully justified since it provides a very good approximation of the effective jump conditions.

\section{Conclusions}
\label{sec:conclusions}

The general plasma model  derived in \cite{Bible} has been presented in a simplified model,
 without considering the electro-magnetic forces and under several assumptions, which inherits the difficulties of the general case in terms of evaluating jump conditions and simulating shock solutions.
 We have proposed a decoupling of the governing equations in order to derive an analytic expression of the missing jump condition on the electron temperature; this decoupling even if seemingly non-physical, will prove to provide a good estimate of the jump conditions for the full system in a non-negligible range of Mach number close to one. We observe that this analytic jump condition is rather different from the jump conditions obtained in the literature and the discrepancies get worse as the Mach number increases. In order to reproduce numerically the structure of the traveling wave solution with the proper jump conditions, we have used a finite volume method of Godunov type. A naive consistant treatment of the nonconservative product proves that, for a fine resolution, the expected wave is well reproduced. We have verified the jump condition as well as the structure of the traveling wave obtained analytically and identified the required level of resolution in order to prevent the appearance of an additional artificial jump due to the numerical dissipation of the numerical scheme.
 
In this context, we have developed a numerical scheme with a specific treatment of the nonconservative product. The idea is to express the compatibility equations at the discontinuity of the traveling wave in a discretized sense. It gives the ability to predict the proper traveling wave even when the gradients are not fully resolved. We have thus built a scheme, which is able to capture the traveling wave in the highly- or weakly-resolved cases.
Such a scheme is important since the weakly-discretized case is particularly relevant in the sun chromosphere conditions. We have also applied a Strang operator splitting  technique in order to prevent the use of small time-steps limited by the presence of large diffusion terms and a Fourier stability condition.
Eventually, a 1D traveling wave test case has been presented based on conditions found in the sun chromosphere, which allowed us to assess the numerical scheme and numerical strategy based on operator splitting. Such a strategy should prove very useful for applications in solar physics in order to gain computational cost, as well as obtain physically accurate simulations.

This scheme has then be validated numerically with a specific choice for the initial conditions. Indeed, we have focused on particular solutions of the systems \eqref{eq:fullprob}, \eqref{eq:couple1}, and \eqref{eq:couple2} as traveling waves.
We have provided a numerical treatment of the nonconservative product that can be used in both cases \caseHD and \caseWD, that is even if the length \(\LD\) is much smaller than \(\LT\).
However, the behavior of this scheme was not investigated for more complicated initial conditions and in particular in the case where several traveling waves interact. 

Furthermore, jump conditions based on our knowledge of the decoupled problem have been derived numerically for the fully coupled system and we were able to justify that studying the decoupled problem provides a good approximation of the original system in a range of Mach number close to one, as well as a good insight on the resolution of the problem for the general system of equations. 

Two natural questions arise from this piece of work. A second order version of our first order scheme is highly desirable for applications. Besides, it would be interesting to derive an asymptotic preserving (AP) scheme, when the mass diffusion of electrons is going to zero. Let us underline that the thermal conductivity of electrons if usually high and the limit of no dissipation for the electronic variable is irrelevant from a physical perspective. however,  our compatibility conditions exhibit the fact that the solution smoothness is changing when the electron mass diffusion coefficient is taken as zero and deriving an AP scheme for this limit would be interesting. The right path to resolve these two problems is to obtain a numerical scheme with the ability to decouple completely the discretized version of the compatibility conditions and the consistency condition of the scheme. We are currently investigating this issue in collaboration with F. Coquel. 
The contributions proposed in this paper should also be extended to the case where electro-magnetic forces are present and the system is coupled to Maxwell's equations. 
Once this is accomplished, we have developed a massively parallel code and implemented a general multicomponent model for partially ionized plasma flows coupled to Maxwell equations in \cite{Wargnier_2018} and the extension to multi-D configurations and more complex flows should be rather straightforward. 
These are the subjects of our current research.

\newpage

\appendix

\section{General governing equations} 
\label{sec:annexe}

Multicomponent nonequilibrium Navier-Stokes equations are obtained from the fluid model derived in \cite{Bible,chimie} for the fully magnetized case and the Maxwellian regime for reactive collisions.  The general governing equations read in non dimensional form:
\begin{equation}
\left\lbrace
\begin{aligned}
&\dt\rhoe + \dx\dscal(\rhoe\vitesse) =  -\dx\dscal(\rhoe\Ve) -\epsilon\bigl(\dx\dscal(\rhoe\Vde)- \omegaeo\bigr), \\
&\dt\rhoi + \dx\dscal(\rhoi\vitesse) = -\epsilon\bigl(\dx\dscal(\rhoi\Vi)- \mi\,\omegaio\bigr),\qquad\qquad i \in \lourd,   \\
& \!\begin{multlined}[b]
\dt(\rhoH\vitesse+\epsilon_0\E\pvect\B)+\dx\dscal \bigl(\rhoH\vitesse\ptens\vitesse+ (\pression+\energie^{EM})\identite-(\epsilon_0\E\ptens\E+\tfrac{1}{\mu_0}\B\ptens\B)\bigr) \\
=-\epsilon\dx\dscal \visqueux[\heavy], 
\end{multlined}\\
&\dt(\rhoe\energiee)+\dx\dscal(\rhoe\energiee\vitesse)+\pre\dx\dscal\vitesse=- \dx\dscal \heate+S^{(1)}-\epsilon\bigl(\dx\dscal\heatde-S^{(2)}\bigr), \\
& \!\begin{multlined}[b]
\dt(\energie+\energie^{EM})+\dx\dscal(\enthalpie\vitesse+\tfrac{1}{\mu_0}\E\pvect\B)=-\dx\dscal\heate \qquad\qquad\qquad\qquad\phantom{ }\\
-\epsilon\bigl(\dx\dscal\heatde + \dx\dscal \heati[\heavy] +\dx\dscal(\visqueux[\heavy]\dscal \vitesse)\bigr),
\end{multlined}
\end{aligned}
\right.
\label{eqfullsys}
\end{equation}
with the same notations as in System~\cref{eq:fullprob}. Some new quantities are introduced: symbol  
$\Vde$ stands for the second-order electron diffusion velocity, 
$\omegaeo$ the electron chemical production rate, 
$\rhoi$ the density of the heavy particle $i \in \lourd$,  
$\Vi$ the  diffusion velocity of the heavy particle $i \in \lourd$, 
$\mi$ the non-dimensional mass of the heavy particle $i \in \lourd$, 
$\omegaio$ the chemical production rate of the heavy particle $i \in \lourd$,  
$\epsilon_0$ the vacuum permittivity,  
$\E$ the electric field, 
$\B$ the magnetic field, 
$\energie^{EM}={\epsilon_0 \E^2}/{2}+\B^2/(2 \mu_0 )$ the electromagnetic energy,  $\mu_0$ the vacuum permeability, 
$\visqueux[\heavy]$ the heavy-particle viscous stress tensor,   
\smash{$\heatde$} the second-order electron heat flux, 
$\heati[\heavy]$ the heavy-particle  heat flux.
The source terms $S^{(1)}$ and $S^{(2)}$ are defined as
\begin{equation}
S^{(1)}=\JJe\dscal\Ep-\deltaEho  ,\quad S^{(2)}=\Jde\dscal\Ep-\deltaEhu -\Delta E^{(2)}_{\textup{chem}}, 
\end{equation} 
where quantity $ \JJe=\ne \qe \Ve$ is the first-order electron conduction current density with $\qe$ the electron charge, $\Ep=\E+\vitesse\pvect\B$, the electric field expressed in the heavy-particle reference frame, \smash{$\deltaEho$},
the energy transferred from heavy particles to electrons at order zero, 
$ \Jde= \ne\qe\Vde$, the second-order electron conduction current density,
\smash{$\deltaEhu$}, the energy transferred from heavy particles to electrons at first order zero, and \smash{$\Delta E^{(2)}_{\textup{chem}}$}, the chemistry-energy coupling term. In solar physics applications,	the full system of equations \eqref{eqfullsys} can be coupled with Maxwell's equations \cite{thompson}	
\begin{equation}
\begin{cases}
\begin{split}
&\dx\dscal\E = \frac{nq}{\epsilon_0}, \\
&\dx\dscal\B=0, \\
&\dt\B+\dx\pvect \E=0, \\
& \mu_0\epsilon_0\dt \E- \dx\pvect \B=-\mu_0 \courantel, \\
\end{split}
\end{cases}
\end{equation}
where quantity	$nq$ is the mixture charge, and $\courantel=n q \vitesse+ \JJe +\Jde+\JJi[\heavy] $,  the total current density with 
$\JJi[\heavy]=\sum_{i\in\lourd}\ni\qi\Vi$, the heavy-particle conduction current density with  $\qi$ the charge of the heavy particle $i \in \lourd$. The electron transport fluxes $\Ve$ and $\heate$ appear at the convective time scale corresponding to the Euler equations for the heavy species (zeroth order). The transport fluxes \smash{$\Vde$ and $\heatde$} are obtained at the dissipative time scale corresponding to the Navier-Stokes equations for the heavy species. Notice that electrons  participate in the momentum balance through the pressure gradient and the Lorentz force but they do not contribute to the viscous stress tensor due to their small mass. The electron transport properties are anisotropic and depend on the direction of the magnetic field whereas the heavy transport properties are isotropic. For example, the first-order electron diffusion velocity is expressed by means of Fick's law and Soret's law
\begin{equation}
\Ve = -\mDee(\de+\mchie\glogTe),
\end{equation}
 where $\mDee$ is the electron diffusion coefficient tensor, and $\mchie$ the thermal diffusion ratio. The electron diffusion driving force is $\de={\dx \pre}/{\pre}-{\ne\qe}\Ep /{\pre} $. The first-order electron heat flux is expressed by means of Fourier's law, together with Dufour's law and a term of diffusion of enthalpy
\begin{equation}
\heate=-\mlambdae\dx \tempe+\pre\mchie\Ve+\rhoe \enthalpiee\Ve.
\end{equation}
Quantity $\mlambdae$ is the electron thermal conductivity tensor. These transport properties can be computed at the microscopic level using a spectral Galerkin method.

\section{Convervative models
\label{sec:annexe2}} 

We introduce two additional models.

\textbf{Model} $\mathcal{M}^{\text{ent}}$ with a conservation equation of entropy:
\begin{equation}
\left\lbrace
\begin{aligned}
&\dt (\rhoH) +\dx (\rhoH \vitesse) = 0,   \\
&\dt (\rhoH \vitesse) + \dx (\rhoH \vitesse^{2} + \pression) = 0, \\
&\dt (\energie) + \dx (\energie \vitesse +\pression \vitesse) = 0, \\
&\dt (\rhoe) + \dx (\rhoe \vitesse) = 0, \\
&\dt (\rhoe \entropiee) +\dx(\rhoe\entropiee\vitesse) =0, 
\end{aligned}
\right.
\label{eq:entropy}
\end{equation}
where the  electron entropy $\entropiee$ is defined by the relation
\(\pre = (\gamma-1){\rhoe}^{\!\!\gamma}\exp(\entropiee/c_v)\),
where $c_v$ is the electron specific  heat at constant volume.

\textbf{Model} $\mathcal{M}^{\text{src}}$, with the nonconservative product as a source term:
\begin{equation}
\left\lbrace
\begin{aligned}
\dt& (\rhoH) +\dx (\rhoH \vitesse) = 0,   \\
\dt& (\rhoH \vitesse) + \dx (\rhoH \vitesse^{2} + \pression) = 0, \\
\dt& (\energie) + \dx (\energie \vitesse +\pression \vitesse) = 0, \\
\dt& (\rhoe) + \dx (\rhoe \vitesse) = 0, \\
\dt& (\rhoe \energiee) +\dx(\rhoe\energiee\vitesse) =0.
\end{aligned}
\right.
\label{eq:source}
\end{equation}

\bibliographystyle{siamplain}
\bibliography{references}
\end{document}